\documentclass[10pt,amstex]{article}
\usepackage{amsmath,amssymb}
\usepackage{enumerate}
\usepackage{amsthm}
\date{\empty}

\begin{document}
\begin{center}
{\bf\huge  2-microlocal spaces associated with   
Besov type and Triebel--Lizorkin type spaces}
\end{center}

\bigskip

\bigskip

\begin{center}
{\bf Koichi Saka}
\end{center}

\bigskip

\bigskip


\bigskip

\bigskip


\bigskip

\bigskip

\noindent
{\bf Abstract}. 
In this paper we introduce and investigate new 2-microlocal spaces associated with Besov type and Triebel--Lizorkin type spaces.  
 We establish characterizations of these  function spaces via 
the $\varphi$--transform,  the atomic and molecular decomposition  
 and the wavelet decomposition. As applications we consider 
boundedness of 
 the Calder$\acute{\rm o}$n--Zygmund operator and the pseudo--differential operator on the function spaces.

\bigskip

\noindent
{\bf Key words}:

\bigskip

Wavelet,  Besov space,  Triebel--Lizorkin space, 
 Calder$\acute{\rm o}$n--Zygmund operator, pseudo--differential operator, $\varphi$-transform, atomic and molecular decompostion,  2-microlocal space.

\bigskip

\noindent
{\bf 2010 Mathematics Subject Classification}:

\bigskip

\begin{tabular}{|c|c|c|c|c|c|}
\hline
4&2&B&3&5 \\
\hline
\end{tabular}(Primary),
     42B25, 42C40 (Secondary).

\newpage

\begin{center}
{\bf 1. Introduction}
\end{center}

It is well known that function spaces have increasing applications in  
many areas of modern analysis, in particular,  harmonic analysis and 
partial differential equations. 
The most general scales, probably, 
 are the scales of Besov spaces and Triebel--Lizorkin spaces which cover many 
 well known classical concrete function spaces such as Lebesgue spaces, 
 Lipschitz spaces, Sobolev spaces, Hardy spaces and  BMO spaces ([23,24]).

 D. Yang and W. Yuan in [25] and W. Sickel, D. Yang and W. Yuan in [22], 
  introduced a class of Besov type and Triebel--Lizorkin type spaces 
  which generalized many classical function spaces such as Besov spaces, Triebel--Lizorkin spaces, Morrey spaces and $Q$-type spaces. 

The 2-microloal spaces are due to Bony [2] in order to study the propagation of singularities of the solutions of nonlinear evolution equations. 
It is an appropriate instrument to describe the local regularity and the oscillatory behavior of functions near to singularity (Meyer [18]). The theory has been elaborated and widely used in fractal analysis and signal processing. For systematic discussions of the concept and further references of 2-microlocal spaces, we refer to Meyer[18,19],  Levy-Vehel and Seuret [17], Jaffard [10,11]  and Jaffard and M${\rm \acute{e}}$lot [12]. 
 
 The 2-microlocal spaces have been generalized by Kempka [13,14] and Jaffard [10,11] in two different directions. Indeed, Kempka generalized as 2-microlocal Besov spaces with  more general  admissible weights ([13,14]) and Jaffard generalized as a general pointwise regularity associated with Banach spaces or quasi-Banach spaces ([9,10,11]). In this paper we introduce new 2-microlocal spaces based 
on Jaffard's idea and   we will investigate the properties and characterizations of 
these new 2-microlocal spaces associated with Besov type and Triebel--Lizorkin type spaces   which contain many classical function spaces such as 
Besov type spaces, Triebel--Lizorkin type spaces and Morrey spaces as examples . 
 We also study 
relations between  the 
new 2-microlocal  spaces and the classical 2-microlocal spaces. 
Moreover, we consider boundedness of 
 the Calder$\acute{\rm o}$n--Zygmund operator and the pseudo--differential operator on the function spaces.

The plan of the remaining sections in the paper is as follows:

In Section 2 we give the definitions of our  new 2-microlocal  spaces  and the notations which are used later and  we give examples for these spaces.
 
In Section 3 we define corresponding sequence spaces 
of our function spaces 

Furthermore, we discuss the almost diagonality and 
we give the conditions under which the almost diagonal operators are bounded on 
the corresponding sequence spaces.
We give  some auxiliary lemmas which are needed for later sections.

In Section 4 we will characterize our function spaces 
through the corresponding sequence spaces
by the atomic and molecular decomposition, the $\varphi$-transform in the sense of Frazier--Jawerth [6] 
 and 
the wavelet decomposition. Moreover, we investigate  properties for these  
function spaces and we observe relations between  the 
new 2-microlocal  spaces and the classical 2-microlocal spaces. 
. 

In Section 5, as applications, we give the conditions under which the Calder$\acute{\rm o}$n--Zygmund operator and the pseudo--differential operator are bounded on the function spaces.

Throughout the paper, we use $C$ to denote a positive constant different in each occasion. But it will depend on the parameters appearing in each assertion. 
The same notation $C$ are 
not necessarily the same on any two occurrences.  
 We use  the notations $i\vee j=\max\{i,\ j\}$, $i\wedge j=\min\{i,\ j\}$ and $s_{+}=\max\{s,\ 0\}$.
The symbol $X \approx Y$ 
 means that there exist positive constants $C_{1}$ and $C_{2}$ such that 
$X \leq C_{1} Y$ and $Y \leq C_{2}X$. 

\bigskip

\begin{center}
{\bf 2. Definitions}
\end{center}

We consider the dyadic cubes in ${\Bbb R}^n$ of the form 
 $Q=[0,\  2^{-l})^n + 2^{-l}k$ for $k \in {\Bbb Z}^n$ and $ l \in {\Bbb Z}$ 
and use  the notations $l(Q)=2^{-l}$ for the side length and $x_{Q}= 2^{-l}k$ 
for the corner point. Throughout the paper, we use  the notations $P,\ Q,\ R\ $for the dyadic cubes of the form $[0,\  2^{-l})^n + 2^{-l}k$ in ${\mathbb R}^n$, and when the dyadic cubes $Q$ appear as indices, it is understood that $Q$ runs over the  dyadic cubes of this form in ${\mathbb R}^n$. 
We denote by $\mathcal{D}$ the set of all dyadic cubes of this form. 
For a dyadic cube $Q$ and a constant $c>1$, $cQ$ denotes the cube of same center as  $Q$ and $c$ times larger.
We denote by $\chi_{E}$ the characteristic function of a set $E$ 
in ${\mathbb R}^n$. 

We set $\mathbb{N}= 
\{ 1, 2, \cdots \}$ and $\mathbb{N}_{0} = \mathbb{N} \cup \{ 0 \}$.
 Let ${\mathcal S}={\mathcal S}({\Bbb R}^n)$ be the space of 
all Schwartz functions on 
${\Bbb R}^n$ and ${\mathcal S}'$ its dual. 
 Let 
${\mathcal S}_{\infty} ={\mathcal S}_{\infty}({\Bbb R}^n)=  
 \{ f \in {\mathcal S}: \int_{\mathbb{R}^n} f(x) x^{\gamma} dx =0 
\ {\rm for}\ {\rm all}\ \gamma \in {\mathbb N}_{0}^{n} \}$, 
${\mathcal S}_{r} ={\mathcal S}_{r}({\Bbb R}^n)=  
 \{ f \in {\mathcal S}: \int_{\mathbb{R}^n} f(x) x^{\gamma} dx =0 
\ {\rm for}\ {\rm all}\ |\gamma| \leq r  \}$,\ $r \in \mathbb{N}_{0}$. 
We recall that  
the dual space 
${\mathcal S}'_{\infty}={\mathcal S}'_{\infty}({\Bbb R}^n)$
  can be identified with the space 
${\mathcal S}'/ \mathfrak{P}$ of equivalence classes of distributions modulo 
the space $\mathfrak{P}$ of all polynomials on ${\Bbb R}^n$, 
and the dual space  ${\mathcal S}'_{r}={\mathcal S}'_{r}({\Bbb R}^n)$
  can be identified  with the space 
${\mathcal S}'/ \mathfrak{P}_{r}$ of equivalence classes of distributions modulo the space $\mathfrak{P}_{r}$ of all polynomials with $\deg \leq r$ on ${\Bbb R}^n$.
 
We use $\langle f,\ g \rangle$ for the standard inner product $\int f \bar{g}$ 
of two functions and the same notation is employed for the action of a distribution $f \in \mathcal{S}'$
on $ \bar{g} \in \mathcal{S}$.

\noindent
Let $\phi$ be a Schwartz function such that 

\bigskip

(1.1)\ \  supp $\hat{\phi} \subset 
\{ \xi \in {\Bbb R}^n: \ \frac{1}{2} \leq |\xi| \leq 2 \}$,

(1.2)\ \ $|\hat{\phi}(\xi)| \geq C > 0$ if $\frac{3}{5} 
\leq |\xi| \leq \frac{5}{3}$. 

\bigskip

\noindent
 We set 
$\phi_{i}(x)= 2^{in}\phi(2^i x)$.  
Let $s \in {\Bbb R}$ and $\phi$ satisfy (1.1) and (1.2) 
as the above.

\noindent
For $f \in {\mathcal S}'_{\infty}$
we define  some sequences 
indexed by dyadic cubes $P$:

\bigskip

$c(\dot{B}^s_{pq})(P)
=(\sum_{i \geq -\log_{2}l(P)}||2^{is}
\phi_{i}*f||_{L^p(P)}^{q})^{1/q}$,\ 
$0 < p,q \leq \infty$, 

$c(\dot{F}^s_{pq})(P)=
||\{\sum_{i \geq -\log_{2}l(P)}(2^{is}
|\phi_{i}*f|)^q\}^{1/q}||_{L^p(P)}$, 
 
$0< p < \infty,\ 0 < q \leq \infty$,

$c(\dot{F}^s_{\infty q})(P)=l(P)^{-\frac{n}{q}}
||\{\sum_{i \geq -\log_{2}l(P)}(2^{is}
|\phi_{i}*f|)^q\}^{1/q}||_{L^q(P)}$,

$0 < q \leq \infty$

\noindent
 with the usual modification for $q=\infty$ respectively.

\noindent
We use the notation  
 $\dot{E}^s_{pq}$ to denote either $\dot{B}^s_{pq}$ or 
$\dot{F}^s_{pq}$.

\bigskip

{\bf Definition 1.}\ \ 
Let $s, \ s' ,\ \sigma\ \in {\Bbb R},\ 0 < p,q \leq \infty$ and  $x_{0} \in \mathbb{R}^n$.

The space $A^s(\dot{E}^{s'}_{pq})^{\sigma}_{x_{0}}$ is defined to be the space 
of all $f \in {\mathcal S}'_{\infty}$ such that 

$$
 ||f||_{A^s(\dot{E}^{s'}_{pq})^{\sigma}_{x_{0}}} 
\equiv
\sup_{\mathcal{D} \ni Q \ni x_{0} }l(Q)^{-\sigma}
\sup_{ \mathcal{D}\ni P \subset 3Q }l(P)^{-s}
c(\dot{E}^{s'}_{pq})(P) < \infty.
$$ 

\noindent
We use the abbreviation 
$A^0(\dot{E}^{s'}_{pq})^{\sigma}_{x_{0}} \equiv 
(\dot{E}^{s'}_{pq})^{\sigma}_{x_{0}}$, 
$A^s(\dot{E}^{s'}_{pq})^{0}_{x_{0}}\equiv A^s(\dot{E}^{s'}_{pq})$ and 
$A^0(\dot{E}^{s'}_{pq})^{0}_{x_{0}}\equiv \dot{E}^{s'}_{pq} 
\equiv \dot{E}^{s'}_{pq}(\mathbb{R}^n)$.  
We note that 
the space $A^s(\dot{E}^{s'}_{pq})$ is the homogeneous Besov type space or the homogeneous Triebel--Lizorkin type space  in the sense of Yang--Yuan[25] and 
the space 
$\dot{E}^{s'}_{pq}\equiv \dot{E}^{s'}_{pq}(\mathbb{R}^n)$ is the classical homogeneous Besov or the classical homogeneous Triebel--Lizorkin space.

\noindent
For 
$f \in {\mathcal S}'_{\infty}$
we define  some sequences 
indexed by dyadic cubes $P$:

\bigskip

$c(\tilde{\dot{B}}^{s'}_{pq})^{\sigma}_{x_{0}}(P)=\\
(\sum_{i \geq -\log_{2}l(P)}||2^{is'}
|\phi_{i}*f(x)|(2^{-i}+|x_{0}-x|)^{-\sigma}||_{L^p(P)}^{q})^{1/q}$,\\ 
$0 < p,\ q \leq \infty$, 

$c(\tilde{\dot{F}}^{s'}_{pq})^{\sigma}_{x_{0}}(P)=\\
||\{\sum_{i  \geq -\log_{2}l(P) }(2^{is'}
|\phi_{i}*f(x)|(2^{-i}+|x_{0}-x|)^{-\sigma})^q\}^{1/q}||_{L^p(P)} $,\\ 
 $0< p < \infty$,\ $0 <  q \leq \infty$,

$c(\tilde{\dot{F}}^{s'}_{\infty q})^{\sigma}_{x_{0}}(P)=\\
l(P)^{-\frac{n}{q}}
||\{\sum_{i \geq -\log_{2}l(P)}(2^{is'}
|\phi_{i}*f(x)|(2^{-i}+|x_{0}-x|)^{-\sigma})^q\}^{1/q}||_{L^q(P)}$,\\
$0 <  q \leq \infty$

\noindent
with the usual modification for $q=\infty$ respectively.

\noindent
We  use the notation $\tilde{\dot{E}}^{s'}_{pq}$ to denote either 
$\tilde{\dot{B}}^s_{pq}$ or 
$\tilde{\dot{F}}^s_{pq}$.

\bigskip

{\bf Definition 2.}\ \ 
Let $s, \ s', \ \sigma\ \in {\Bbb R},\  0< p,q \leq \infty$ and 
  $x_{0} \in \mathbb{R}^n$.
 
The space $A^s(\tilde{\dot{E}}^{s'}_{pq})^{\sigma}_{x_{0}}$ is defined to be the space of  all $f \in \mathcal{S}_{\infty}'$ such that 
$$
||f||_{A^s(\tilde{\dot{E}}^{s'}_{pq})^{\sigma}_{x_{0}}}\equiv \sup_{\mathcal{D}\ni P}l(P)^{-s} c(\tilde{\dot{E}}^{s'}_{pq})^{\sigma}_{x_{0}}(P) < \infty.
$$

The space $A^s(\tilde{\dot{E}}^{s'}_{pq})^{\sigma}_{x_{0}}$ is the classical
2 microlocal Besov or the classical 2 microlocal Triebel--Lizorkin space.

\noindent
We use the abbreviation 
$A^0(\tilde{\dot{E}}^{s'}_{pq})^{\sigma}_{x_{0}} \equiv 
(\tilde{\dot{E}}^{s'}_{pq})^{\sigma}_{x_{0}}$. 

\bigskip

{\bf Examples.}

\begin{enumerate}[(i)]
\item \  The spaces $A^0(\dot{E}^{s'}_{pq})^{0}_{x_{0}}=A^0(\tilde{\dot{E}}^{s'}_{pq})^{0}_{x_{0}}= \dot{E}^{s'}_{pq}(\mathbb{R}^n)$ 
 are the  homogeneous Besov spaces or the 
homogeneous Triebel--Lizorkin spaces,
[23].

\item \  
The Besov type spaces $\dot{B}^{s, \tau}_{pq}(\mathbb{R}^n)$  and 
the Triebel--Lizorkin type spaces $\dot{F}^{s, \tau}_{pq}(\mathbb{R}^n)$ 
introduced by D. Yang and W. Yuan [25]  ,
are contained in our definition as special cases that

\noindent
$\dot{E}^{s', s}_{pq}(\mathbb{R}^n)=A^{ns}(\dot{E}^{s'}_{pq})^{0}_{x_{0}}=
A^{ns}(\tilde{\dot{E}}^{s'}_{pq})^{0}_{x_{0}}
$.

\item \  The Besov-Morrey  $\dot{{\mathcal N}}^s_{uqp}$ and the Triebel--Lizorkin-Morrey spaces   $\dot{{\mathcal E}}^s_{uqp}$ 
studied by  Y. Sawano and H. Tanaka [20], or Y. Sawano, D. Yang and W. Yuan [21] are realized in our definition as  

$\dot{{\mathcal N}}^s_{uqp} \subset A^{n(\frac{1}{p}-\frac{1}{u})}
(\dot{B}^s_{pq})^0_{x_{0}}$ if $0 < p \leq u \leq \infty$ and $0 < q \leq \infty$,

$\dot{{\mathcal E}}^s_{uqp}=A^{n(\frac{1}{p}-\frac{1}{u})}
(\dot{F}^s_{pq})^0_{x_{0}}$ if $0 < p \leq u \leq \infty$ and $0 < q \leq \infty$, and

Especially the Morrey space ${\mathcal M}^u_{p}$ is realized as

${\mathcal M}^u_{p}= 
A^{n(\frac{1}{p}-\frac{1}{u})}(\dot{F}^0_{p2})^0_{x_{0}}$,
$1 < p < u < \infty$.

\item \  The $\dot{B}_{\sigma}$-Morrey spaces $\dot{B}_{\sigma}(L_{p,\lambda})$
studied by Y. Komori-Furuya et al. [15], 
are contained in our definition as special cases, 
that is, 
 
$\dot{B}_{\sigma}(L_{p,\lambda})$
=$A^{\lambda+\frac{n}{p}}(\dot{F}^0_{p2})^{\sigma}_{0}$, \ $1 < p < \infty$.

\item \ The local Morrey spaces   
 $LM_{p, \lambda}$ 
studied by Ts. Batbold and Y. Sawano [1], 
are realized in our cases as 

$LM_{p, \lambda}=(\dot{F}^{0}_{p2})^{\lambda/p}_{0}$
, \ $1 < p < \infty$.

\end{enumerate}

\bigskip

\begin{center}
{\bf 3. Sequence spaces}
\end{center}

For a sequence $c=(c(P))$ we define  some sequences indexed by dyadic cubes 
$P$:

\bigskip

$
c(\dot{b}^s_{pq})(P)
=(\sum_{i \geq -\log_{2}l(P)}||\sum_{l(R)=2^{-i}}2^{is}
|c(R)|\chi_{R}||_{L^p(P)}^{q})^{1/q}$, 

$ 0 < p,\ q \leq \infty$,

$
c(\dot{f}^s_{pq})(P)=||\bigl\{\sum_{i \geq -\log_{2}l(P)}\bigl(
\sum_{l(R)=2^{-i}}2^{is}
|c(R)|\chi_{R}\bigr)^q\bigr\}^{1/q}||_{L^p(P)}$,

$ 0 < p < \infty,\ \ 0 < q \leq \infty$,
 and

$
c(\dot{f}^s_{\infty q})(P)=l(P)^{-\frac{n}{q}}\times
\\
||\bigl\{\sum_{i \geq -\log_{2}l(P)}\bigl(\sum_{l(R)=2^{-i}}2^{is}
|c(R)|\chi_{R}\bigr)^q\bigr\}^{1/q}||_{L^q(P)}$,
$ 0 < q \leq \infty$, 

\noindent
with the usual modification for $q=\infty$ respectively.
 
\noindent
We  use the notation 
 $\dot{e}^s_{pq}$ to denote either $\dot{b}^s_{pq}$ or 
$\dot{f}^s_{pq}$.

\bigskip

{\bf Definition 3.}\ \ 
Let $s, \ s' ,\ \sigma\ \in {\Bbb R},\ 0< p,q \leq \infty$ and 
 $x_{0} \in \mathbb{R}^n$.
 
 We define the sequence space $a^s(\dot{e}^{s'}_{pq})^{\sigma}_{x_{0}}$ 
to be the space of all sequences $c=(c(P))$ such that 
$$
||c||_{a^s(\dot{e}^{s'}_{pq})^{\sigma}_{x_{0}}}\equiv 
\sup_{\mathcal{D} \ni Q \ni x_{0}}l(Q)^{-\sigma}
\sup_{\mathcal{D} \ni P \subset 3Q}l(P)^{-s}c(\dot{e}^{s'}_{pq})(P)
< \infty.
$$

\noindent
We use the abbreviation 
$a^0(\dot{e}^{s'}_{pq})^{\sigma}_{x_{0}} \equiv 
(\dot{e}^{s'}_{pq})^{\sigma}_{x_{0}}$, $a^s(\dot{e}^{s'}_{pq})^{0}_{x_{0}}\equiv a^s(\dot{e}^{s'}_{pq})$ and 

\noindent
$a^0(\dot{e}^{s'}_{pq})^{0}_{x_{0}}
\equiv \dot{e}^{s'}_{pq}
\equiv 
\dot{e}^{s'}_{pq}(\mathbb{R}^n)$. 
We note that the space $a^s(\dot{e}^{s'}_{pq})$ 
is the sequence space of the homogeneous Besov type space or the homogeneous Triebel--Lizorkin type space  in the sense of Yang--Yuan[25] and 
the space 
$\dot{e}^{s'}_{pq}\equiv \dot{e}^{s'}_{pq}(\mathbb{R}^n)$ 
is the sequence space of the classical homogeneous Besov or the classical homogeneous Triebel--Lizorkin space.

We define a sequential version for Definition 2.

\noindent
We define that for a sequence $(c(R))$,

\bigskip

$c(\tilde{\dot{b}}^{s'}_{pq})^{\sigma}_{x_{0}}(P)=\\
(\sum_{i \geq -\log_{2}l(P)}||
\sum_{l(R)=2^{-i}}2^{is'}|c(R)|(2^{-i}+|x_{0}-x|)^{-\sigma}\chi_{R}||_{L^p(P)}^{q})^{1/q} $, 
\\
$0 < p,\ q \leq \infty$,

$c(\tilde{\dot{f}}^{s'}_{pq})^{\sigma}_{x_{0}}(P)=\\
||\{\sum_{i \geq -\log_{2}l(P)}(\sum_{l(R)=2^{-i}}2^{is'}
|c(R)|(2^{-i}+|x_{0}-x|)^{-\sigma}\chi_{R})^q\}^{1/q}||_{L^p(P)}$, 
\\
 $0< p < \infty$,\ $0 <  q \leq \infty$,

$c(\tilde{\dot{f}}^{s'}_{\infty q})^{\sigma}_{x_{0}}=
l(P)^{-\frac{n}{q}}\times \\
||\{\sum_{i\geq -\log_{2}l(P)}(\sum_{l(R)=2^{-i}}2^{is'}|c(R)|
(2^{-i}+|x_{0}-x|)^{-\sigma}\chi_{R})^q\}^{1/q}||_{L^q(P)}$,
\\
$0 <  q \leq \infty$,

 \noindent
 with the usual modification for $q=\infty$ respectively.

\bigskip

{\bf Definition 4.}
\ \ Let $s, \ s' ,\ \sigma\ \in {\Bbb R},\ 0< p,q \leq \infty$ and 
 $x_{0} \in \mathbb{R}^n$.

We define the sequence space $a^s(\tilde{\dot{e}}^{s'}_{pq})^{\sigma}_{x_{0}}$
to be the space of all sequences c=(c(R)) such that 
$$
||c||_{a^s(\tilde{\dot{e}}^{s'}_{pq})^{\sigma}_{x_{0}}}\equiv 
 \sup_{\mathcal{D} \ni P}l(P)^{-s}c(\tilde{\dot{e}}^{s'}_{pq})^{\sigma}_{x_{0}}(P) < \infty.  
$$

\bigskip

\noindent
We use the abbreviation 
$a^0(\tilde{\dot{e}}^{s'}_{pq})^{\sigma}_{x_{0}} \equiv 
(\tilde{\dot{e}}^{s'}_{pq})^{\sigma}_{x_{0}}$. 

\bigskip

We will give some lemmas which are needed through the whole sections later.

\bigskip

{\bf Definition 5.}\ \ 
Let $r_{1},\ r_{2}  \geq 0$  and $L > 0$. 
We say that a matrix operator $A=\{ a_{QP} \}_{Q P}$, 
indexed by dyadic cubes $Q$ and $P$, 
is  ($r_{1}, r_{2}, L$)-almost diagonal if the matrix $\{ a_{QP} \}$ 
satisfies that 

\bigskip

$|a_{QP}| \leq C\bigl( \frac{l(Q)}{l(P)} \bigl)^{r_{1}}(1+l(P)^{-1}|x_{Q}-x_{P}|)^{-L}$ if $l(Q) \leq l(P)$,

 $
|a_{Q P}| \leq C\bigl( \frac{l(P)}{l(Q)} \bigl)^{r_{2}}(1+l(Q)^{-1}|x_{Q}-x_{P}|)^{-L}$ if $l(Q) > l(P)$.

\bigskip

The results about the boundedness of almost diagonal operators in 
 [6: Theorem 3.3], also
hold in our cases.
  
\bigskip

{\bf Lemma\ 1.}\ \ {\it 
Suppose that $ s,\ s',\ \sigma \in {\Bbb R},\ x_{0} \in \mathbb{R}^n$ 
and 
\ $0 < p,\ q\ \leq \infty$. 
 Then

(i)\ \ an $(r_{1}, r_{2}, L)$-almost diagonal 
matrix operator A is bounded on $a^s(\dot{e}^{s'}_{pq})^{\sigma}_{x_{0}}$ 
for $r_{1} > \max(s',\ \sigma+s + s'- \frac{n}{p}),\ r_{2} > J-s'$ and $L>J$ 
where 
$J= n/ \min(1,\ p, \ q)$ 
in the case $\dot{e}^{s'}_{pq}=\dot{f}^{s'}_{pq}$
 and 
$J= n/ \min(1,\ p)$ in the case $\dot{e}^{s'}_{pq}=\dot{b}^{s'}_{pq}$ 
 respectively.

(ii)\ \ An $(r_{1}, r_{2}, L)$-almost diagonal 
matrix operator A  is bounded on 
$a^s(\tilde{\dot{e}}^{s'}_{pq})^{\sigma}_{x_{0}}$ 
for $r_{1} > \max(s'+(\sigma\vee 0),\  (\sigma\vee 0)+s + s'- \frac{n}{p}),\ r_{2} > J-s'-(\sigma\wedge 0)$ and $L>J$ 
where $J= n/ \min(1,\ p, \ q)$ 
in the case $\tilde{\dot{e}}^{s'}_{pq}=\tilde{\dot{f}}^{s'}_{pq}$
 and 
$J= n/ \min(1,\ p)$ in the case $\tilde{\dot{e}}^{s'}_{pq}=\tilde{\dot{b}}^{s'}_{pq}$ 
respectively.
}
  
\bigskip

{\it Proof} :\ \ (i)\ \ We may assume $\sigma \geq 0$ since 
$a^s(\dot{e}^{s'}_{pq})^{\sigma}_{x_{0}}=\{0\}$ if $\sigma < 0$ (See Proposition 1 below). 
We assume that $A = (a_{R R'})$ is $(r_{1}, r_{2}, L)$- almost diagonal. 
Let $c= (c(R)) \in a^s(\dot{e}^{s'}_{pq})^{\sigma}_{x_{0}}$. 
For  dyadic cubes $P$ and $R$ with $R \subset P$, 
we write $Ac= A_{0}c + A_{1}c + A_{2}c$ with 

\begin{eqnarray*}
&&
(A_{0}c)(R)=\sum_{l(R) \leq l(R') \leq l(P)}a_{R R'}c(R'), 
\\
&&
(A_{1}c)(R)=\sum_{l(R')  < l(R) \leq l(P)}a_{R R'}c(R'), 
\\
&&
(A_{2}c)(R)=\sum_{l(R) \leq l(P) < l(R') }a_{R R'}c(R').
\end{eqnarray*}
We claim that
$$
||A_{i}c||_{a^{s}(\dot{e}^{s'}_{pq})^{\sigma}_{x_{0}}}
\leq 
C||c||_{a^{s}(\dot{e}^{s'}_{pq})^{\sigma}_{x_{0}}}, \ \ i=0,1,2.
$$
First we will consider the case of the F-type for $0 < p < \infty$, $0 < q \leq \infty$.
Since $A$ is almost diagonal, we see that 
for dyadic cubes $P$
 with $l(P)=2^{-j}$, 
\begin{eqnarray*}
\lefteqn{(A_{0}c)(\dot{f}^{s'}_{pq})(P) = 
||\bigl\{\sum_{i \geq j}\sum_{l(R)=2^{-i}}\bigl(2^{is'}|(A_{0}c)(R)|
\bigr)^{q}\chi_{R}\bigr\}^{1/q}||_{L^{p}(P)}
} 
\\
&\leq &
C||\bigl\{\sum_{i \geq j}\sum_{l(R)=2^{-i}}2^{is'q}\bigl(
\sum_{i \geq k \geq j}\sum_{l(R')=2^{-k}}|a_{R R'}||c(R')|\bigr)^{q}
\chi_{R}\bigr\}^{1/q}||_{L^{p}(P)} 
\\
&\leq&
C||\bigl\{\sum_{i \geq j}\sum_{l(R)=2^{-i}}2^{is'q}\times 
\\
\lefteqn{\bigl(\sum_{i \geq k \geq j}
\sum_{l(R')=2^{-k}}2^{-(i-k)r_{1}}(1+2^k|x_{R}-x_{R'}|)^{-L}
|c(R')|\bigr)^{q}\chi_{R}\bigr\}^{1/q}||_{L^{p}(P)}.
}
\end{eqnarray*}
Using the maximal function $M_{t}f(x)$,\ $0< t\leq 1$,  defined by 
$$
M_{t}f(x)=\sup_{x \in Q}(\frac{1}{l(Q)^n}\int_{Q}
|f(y)|^t\ dy)^{1/t}
$$
(cf. [16: Lemma 7.1] or [6: Remark A.3]), 
we have for $L > n/t$,
\begin{eqnarray*}
\lefteqn{
(A_{0}c)(\dot{f}^{s'}_{pq})(P) \leq 
C||\bigl\{\sum_{i \geq j}\sum_{l(R)=2^{-i}}2^{is'q}2^{-ir_{1}q}\times
}
\\
&&
\Bigl(\sum_{i \geq k \geq j}2^{k r_{1}}2^{(k-i)_{+}n/t}M_{t}\bigl(\sum_{l(R')=2^{-k}}|c(R')|
\chi_{R'}\bigr)
\Bigr)^{q}\chi_{R}\bigr\}^{1/q}||_{L^{p}(P)} 
\\
&\leq& 
C||\big\{\sum_{i \geq j}2^{-i(r_{1}-s')q}\Bigl(\sum_{i \geq k \geq j}
2^{kr_{1}}
M_{t}\bigl(\sum_{l(R')=2^{-k}}|c(R')|\chi_{R'}\bigr)\Bigr)
^{q}\bigr\}^{1/q}||
_{L^{p}(P)} 
\\
&\leq& 
C||\bigl\{\sum_{i \geq j}2^{is'q}M_{t}\bigl(\sum_{l(R')=2^{-i}}|c(R')|
\chi_{R'}\bigr)^{q}
\bigr\}^{1/q}||_{L^{p}(P)} 
\\
&\leq &
C||\big\{\sum_{i \geq j}2^{is'q}\bigl(\sum_{l(R')=2^{-i}}|c(R')|\chi_{R'}
\bigr)^{q}
\bigr\}^{1/q}||_{L^{p}(P)} 
=Cc(\dot{f}^{s'}_{pq})(P)
\end{eqnarray*}
where these inequalities follow from Hardy's inequality if $r_{1} > s'$ and the  Fefferman-Stein vector valued inequality if $0< t< \min(p,q)$.

For the B-type case  we have for $L > n/t$,
\begin{eqnarray*}
\lefteqn{
(A_{0}c)(\dot{b}^{s'}_{pq})(P) = 
\bigl(\sum_{i \geq j}||\sum_{l(R)=2^{-i}}2^{is'}|(A_{0}c)(R)|
\chi_{R}||^{q}_{L^{p}(P)}\bigr)^{1/q}
}
\\
\lefteqn{ 
\leq C\bigl\{\sum_{i \geq j}||\sum_{l(R)=2^{-i}}2^{is'}\times
}
\\
&&
\sum_{i \geq k \geq j}
\sum_{l(R')=2^{-k}}2^{-(i-k)r_{1}}(1+2^{k}|x_{R}-x_{R'}|)^{-L}
|c(R')|\chi_{R}||^{q}_{L^{p}(P)}\bigr\}^{1/q} 
\\
\lefteqn{\leq 
C\bigl\{\sum_{i \geq j}2^{-i(r_{1}-s')q}\times
}
\\
&&
||\sum_{l(R)=2^{-i}}
\sum_{i \geq k \geq j}2^{kr_{1}}M_{t}(\sum_{l(R')=2^{-k}}|c(R')|\chi_{R'})
\chi_{R}||^{q}_{L^{p}(P)}\bigr\}^{1/q} 
\\
\lefteqn{\leq 
C\bigl\{\sum_{i \geq j}2^{-i(r_{1}-s')q}\bigl(\sum_{i \geq k \geq j}2^{kr_{1}}
||M_{t}(\sum_{l(R')=2^{-k}}|c(R')|\chi_{R'})||_{L^{p}(P)}\bigr)^{q}\bigr\}
^{1/q} 
}
\\
\lefteqn{\leq 
C\bigl\{\sum_{i \geq j}2^{-i(r_{1}-s')q}\bigl(\sum_{i \geq k \geq j}2^{kr_{1}}
||\sum_{l(R')=2^{-k}}|c(R')|\chi_{R'}||_{L^{p}(P)}\bigr)^{q}\bigr\}^{1/q} 
}
\\
\lefteqn{\leq
C\bigl(\sum_{i \geq j}2^{is'q}||\sum_{l(R')=2^{-i}}|c(R')|\chi_{R'}||^{q}
_{L^{p}(P)}\bigr)^{1/q} 
=Cc(\dot{b}^{s'}_{pq})(P)
}
\end{eqnarray*}

\noindent
where the inequalities follow from Hardy's inequality if $r_{1} > s'$ and the  Fefferman-Stein vector valued inequality if $0 < t< \min(1,p)$.

Hence we get the estimate 

$$
A_{0}c(\dot{e}^{s'}_{pq})(P)
\leq 
Cc(\dot{e}^{s'}_{pq})(P)
$$
if $r_{1} > s'$,  $0 < p < \infty$, $0 < q \leq \infty$, $L > J$.

Similarly  we will see the estimate for $(A_{1}c)(\dot{f}^{s'}_{pq})(P)$.
 We have that for dyadic cubes $P$ with $l(P)=2^{-j}$, 
\begin{eqnarray*}
\lefteqn{(A_{1}c)(\dot{f}^{s'}_{pq})(P) = 
||\bigl\{\sum_{i \geq j}\sum_{l(R)=2^{-i}}\bigl(2^{is'}|(A_{1}c)(R)|
\bigr)^{q}\chi_{R}\bigr\}^{1/q}||_{L^{p}(P)}
} 
\\
&\leq &
C||\bigl\{\sum_{i \geq j}\sum_{l(R)=2^{-i}}2^{is'q}\bigl(
\sum_{i \leq k }\sum_{l(R')=2^{-k}}|a_{R R'}||c(R')|\bigr)^{q}
\chi_{R}\bigr\}^{1/q}||_{L^{p}(P)} 
\\
&\leq&
C||\bigl\{\sum_{i \geq j}\sum_{l(R)=2^{-i}}2^{is'q}\times 
\\
\lefteqn{\bigl(\sum_{i \leq k }
\sum_{l(R')=2^{-k}}2^{-(k-i)r_{2}}(1+2^i|x_{R}-x_{R'}|)^{-L}
|c(R')|\bigr)^{q}\chi_{R}\bigr\}^{1/q}||_{L^{p}(P)}.
}
\end{eqnarray*}
Using the maximal function $M_{t}f(x)$ as the above, 
we have 
\begin{eqnarray*}
\lefteqn{
(A_{1}c)(\dot{f}^{s'}_{pq})(P) \leq 
C||\bigl\{\sum_{i \geq j}\sum_{l(R)=2^{-i}}2^{is'q}2^{ir_{2}q}\times
}
\\
&&
\Bigl(\sum_{i \leq k }2^{-k r_{2}}2^{(k-i)_{+}n/t}M_{t}\bigl(\sum_{l(R')=2^{-k}}|c(R')|
\chi_{R'}\bigr)
\Bigr)^{q}\chi_{R}\bigr\}^{1/q}||_{L^{p}(P)} 
\\
&\leq& 
C||\big\{\sum_{i \geq j}2^{i(r_{2}+s'-n/t)q}\times
\\
&&
\Bigl(\sum_{i \leq k }
2^{-k(r_{2}-n/t)}M_{t}\bigl(\sum_{l(R')=2^{-k}}|c(R')|\chi_{R'}\bigr)\Bigr)
^{q}\bigr\}^{1/q}||
_{L^{p}(P)} 
\\
&\leq& 
C||\bigl\{\sum_{i \geq j}2^{is'q}M_{t}\bigl(\sum_{l(R')=2^{-i}}|c(R')|
\chi_{R'}\bigr)^{q}
\bigr\}^{1/q}||_{L^{p}(P)} 
\\
&\leq& 
C||\big\{\sum_{i \geq j}2^{is'q}\bigl(\sum_{l(R')=2^{-i}}|c(R')|\chi_{R'}
\bigr)^{q}
\bigr\}^{1/q}||_{L^{p}(P)} 
=Cc(\dot{f}^{s'}_{pq})(P)
\end{eqnarray*}
where these inequalities follow from Hardy's inequality if $r_{2} + s'-n/t > 0$
and the  Fefferman-Stein vector valued inequality if $0 < t < \min(p,q)$.

Similarly we get for B-type case 
 that 
$$
(A_{1}c)(\dot{b}^{s'}_{pq})(P) \leq Cc(\dot{b}^{s'}_{pq})(P)
$$
if $r_{2} +s'- n/t>0$, $0 < t < \min(1, p)$. 
Hence we get the estimate 

$$
A_{1}c(\dot{e}^{s'}_{pq})(P)
\leq 
Cc(\dot{e}^{s'}_{pq})(P)
$$
if $r_{2} >J- s'$, 
 $0 < p < \infty$, $0 < q \leq \infty$, $L > J$.

When $p=\infty$, we get the same estimate.
Thus we get 

$$
||A_{i}c||_{a^{s}(\dot{e}^{s'}_{pq})^{\sigma}_{x_{0}}}
\leq 
C||c||_{a^{s}(\dot{e}^{s'}_{pq})^{\sigma}_{x_{0}}}, \ \ i=0,1
$$

\noindent
if $r_{1} > s'$, $r_{2} >J-s'$, $L>J$, $0 < p \leq \infty$ and $0 < q \leq \infty$.
 
Next, we will give the estimates for the $A_{2}$ case.

We note that  if $L > n$ , 
$$\sum_{l(P)=2^{-j}}(1+2^{j}|x_{R}-x_{P}|)^{-L} < \infty$$ 
(cf. [3, Lemma 3.4]), and 

\noindent
for $c \in a^{s}(\dot{e}^{s'}_{pq})^{\sigma}_{x_{0}}$, we have
 $|c(P)| \leq C(|x_{0}-x_{P}|+l(P))^{\sigma}l(P)^{s+s'-n/p}
||c||_{a^{s}(\dot{e}^{s'}_{pq})^{\sigma}_{x_{0}}}$ 
for a dyadic cube $P$. Then 
we obtain, for $0 < p <\infty$, $0 < q \leq \infty$ and a dyadic cube $P$ 
with $l(P)=2^{-j}$, 
\begin{eqnarray*}
\lefteqn{
(A_{2}c)(\dot{f}^{s'}_{pq})(P) =
||\bigl\{\sum_{i \geq j}\sum_{l(R)=2^{-i}}\bigl(2^{is'}|(A_{2}c)(R)|
\bigr)^{q}
\chi_{R}\bigr\}^{1/q}||_{L^{p}(P)} 
}
\\
\lefteqn{
\leq 
C||\bigl\{\sum_{i \geq j}\sum_{l(R)=2^{-i}}2^{is'q}\times
}
\\
\lefteqn{
\bigl(\sum_{j \geq k}\sum_{l(R')=2^{-k}}2^{-(i-k)r_{1}}(1+2^{k}|x_{R}-x_{R'}|
\bigr)^{-L}|c(R')|)^{q}\chi_{R}\bigr\}^{1/q}||_{L^{p}(P)} 
}
\\
\lefteqn{\leq
C2^{-j(r_{1}-s'+n/p)}
\sum_{j \geq k}2^{k(r_{1}-s-s'+\frac{n}{p})}(2^{-k}+|x_{0}-x_{P}|)^{\sigma}
||c||_{a^{s}(\dot{f}^{s'}_{pq})^{\sigma}_{x_{0}}}
}
\\
&\leq& C2^{-j(r_{1}-s'+n/p)}
2^{j(r_{1}-\sigma-s-s'+\frac{n}{p})}(1+2^k|x_{0}-x_{P}|)^{\sigma}||c||_{a^{s}(\dot{f}^{s'}_{pq})^{\sigma}_{x_{0}}}
\\
&\leq& C(2^{-j}+|x_{0}-x_{P}|)^{\sigma}2^{-js}||c||_{a^{s}(\dot{f}^{s'}_{pq})^{\sigma}_{x_{0}}}
\end{eqnarray*}
where these inequalities follow 
if $r_{1} >\sigma +s+s'-\frac{n}{p}$,  $r_{1} > s'$, $ L > n$.

Similarly, for the B-type case we have the same estimate.
\noindent
Hence we have,  
 
$$
||A_{2}c||_{a^{s}(\dot{e}^{s'}_{pq})^{\sigma}_{x_{0}}} 
\leq
C||c||_{a^{s}(\dot{e}^{s'}_{pq})^{\sigma}_{x_{0}}}
$$

\noindent
if $r_{1} > \sigma + s+s'-n/p$,\ $r_{1} > s'$,\ $L>J$, $0 < p < \infty$ and $0 < q \leq \infty$.  

We get the same estimate for the case $p = \infty$. 
Thus we obtain the desired conclusion.

(ii)\ \ We put $w_{i} =(2^{-i}+|x_{0}-x|)^{-\sigma}$. 
We see that  $w_{i}\leq 2^{(i-k)_{+}\sigma}w_{k}$ if  $0 \leq \sigma$,
and 
 $w_{i}\leq 2^{(k-i)_{+}\sigma}w_{k}$ if  $0 > \sigma$. Using these inequalities  
we can prove  the desired result by the same way as the above (i).
\qed

\bigskip

{\bf Lemma 2.}\ \ {\it 
 Let $r_{1}, r_{2} \in \mathbb{N}_{0}, \ L >n $ and $L_{1} > n+r_{1}, L_{2} > n + r_{2}$. 

Assume that  for dyadic cubes $P$ and $R$, 
$\phi_{P}$ and $\varphi_{R}$ 
 are functions on $\mathbb{R}^{n}$ satisfying following properties 
such that 
\begin{eqnarray*}
&&
(2.1) \ \ \int_{\mathbb{R}^n}\phi_{P}(x)x^{\gamma} dx =0\ \  
for\  \ |\gamma| < r_{1},
\\
&&
(2.2)\ \ |\phi_{P}(x)| \leq C(1+l(P)^{-1}|x-x_{P}|)^{-\max(L, L_{1})},
\\
&&
(2.3)\ \ |\partial^{\gamma}\phi_{P}(x)| \leq Cl(P)^{-|\gamma|}(1+l(P)^{-1}|x-x_{P}|)^{-L}\ \ 
\\
&&
for\ \ 0<|\gamma| \leq r_{2},
\\
&&
(2.4)\ \ \int_{\mathbb{R}^n}\varphi_{R}(x) x^{\gamma} dx =0\ \ 
for\ \ |\gamma|< r_{2}, 
\\
&&
(2.5)\ \ |\varphi_{R}(x)| \leq C (1+l(R)^{-1}|x-x_{R}|)
^{-\max(L, L_{2})},
\\
&&
(2.6)\ \ |\partial^{\gamma}\varphi_{R}(x)| \leq C l(R)^{-|\gamma|}
(1+l(R)^{-1}|x-x_{R}|)^{-L} \ \ 
\\
&&for\ \ 0<|\gamma| \leq r_{1}
\end{eqnarray*}
where  (2.1) and  (2.6) are void when $r_{1}=0$, and 
(2.3) and (2.4) are void when  $r_{2}=0$. 
Then we have that

\bigskip

$ l(P)^{-n}|
\langle \phi_{P}\ ,\ \varphi_{R} \rangle| 
\leq C \bigl( \frac{l(P)}{l(R)} \bigl)^{r_{1}}(1+l(R)^{-1}|x_{P}-x_{R}|)^{-L}$

  if \ \ $l(P) \leq l(R)$ and 

 $ l(R)^{-n}|\langle \phi_{P}\ ,\ \varphi_{R} \rangle| 
\leq C \bigl( \frac{l(R)}{l(P)} \bigl)^{r_{2}}(1+l(P)^{-1}|x_{P}-x_{R}|)^{-L}$

if\ \  $l(R) < l(P)$.
}
\bigskip

{\it Proof} :\ \  We refer to [6: Corollary  B.3] , [4: Lemma 6.3] or  [16: Lemma 3.1].
\qed

\bigskip
  
{\bf Lemma\ 3.}\ \ {\it 
Suppose that $ s,\ s',\ \sigma \in {\Bbb R}$,\ $ x_{0} \in \mathbb{R}^n$ 
and 
\ $0 < p,\ q\ \leq \infty$. Let $r_{1}, r_{2} \in \mathbb{N}_{0}$ and $L >n$. 
Assume that functions $\phi_{P},\ \varphi_{R}$ satisfy 
the properties (2.1) through (2.6) as in Lemma 2. Let $J$ as in Lemma 1. 
Then we have

${\rm (i)}$\ \ for a dyadic cube $R$ and a sequence $c \in 
a^s(\dot{e}^{s'}_{pq})^{\sigma}_{x_{0}}$, 

\noindent 
$\sum_{\mathcal{D} \ni P}c(P)\langle \phi_{P}\ ,\ \varphi_{R} \rangle$
 is convergent 
if $r_{1} > J-n-s'$,\ $r_{2}>\sigma+s + s'- \frac{n}{p}$ and $L>\max(J,\ n+\sigma)$,

${\rm (ii)}$\ \ for a dyadic cube $R$ and a sequence $c \in a^s(\tilde{\dot{e}}^{s'}_{pq})^{\sigma}_{x_{0}}$, 

\noindent
 $\sum_{\mathcal{D} \ni P}c(P)\langle \phi_{P}\ ,\ \varphi_{R} \rangle$ is convergent if $r_{1} > J-n-s'-(\sigma\wedge 0)$,\ $r_{2}>\sigma+s + s'- \frac{n}{p}$ and $L>J+\sigma$.

}

\bigskip

{\it Proof} :\ \ 
(i)\ \ we write \ 
$\sum_{\mathcal{D} \ni P}c(P)\langle \phi_{P}\ ,\ \varphi_{R} \rangle
=  I_{0} + I_{1}$ with 

\begin{eqnarray*}
&&
I_{0}=\sum_{l(R) < l(P)}c(P)\langle \phi_{P}\ ,\ \varphi_{R} \rangle, 
\\
&&
I_{1}=\sum_{l(P)  \leq l(R) }c(P)\langle \phi_{P}\ ,\ \varphi_{R} \rangle. 
\end{eqnarray*}
We claim that
$
I_{i}< \infty, \ \ i=0,1.
$

We note that  if $L > n$ , 
$$\sum_{l(P)=2^{-j}}(1+2^{j}|x_{R}-x_{P}|)^{-L} < \infty$$ 
(cf. [3, Lemma 3.4]), and 
$$
|c(P)| \leq C(|x_{P}-x_{0}|+l(P))^{\sigma}l(P)^{s+s'-n/p}
||c||_{a^{s}(\dot{e}^{s'}_{pq})^{\sigma}_{x_{0}}}$$ 
for $c \in a^{s}(\dot{e}^{s'}_{pq})^{\sigma}_{x_{0}}$. 
For a dyadic cube $R$ with $l(R)=2^{-i}$, we have, by Lemma 2
\begin{eqnarray*}
\lefteqn{
|I_{0}| \leq C
\sum_{i > j}\sum_{l(P)=2^{-j}}|c(P)||\langle \phi_{P}\ ,\ \varphi_{R} \rangle|
}
\\
&\leq& 
C\sum_{i > j}\sum_{l(P)=2^{-j}}|c(P)
|2^{-in}2^{(j-i)r_{2}}(1+2^{j}|x_{R}-x_{P}|)^{-L}
\\
&\leq&
C\sum_{i > j}2^{-i(r_{2}+n)}
2^{j(r_{2}-s-s'+\frac{n}{p})}
||c||_{a^{s}(\dot{e}^{s'}_{pq})^{\sigma}_{x_{0}}}\times 
\\
&&
\sum_{l(P)=2^{-j}}(|x_{P}-x_{0}|+2^{-j})^{\sigma}(1+2^{j}|x_{R}-x_{P}|)^{-L}
\\
&\leq&
C\sum_{i > j}2^{-i(r_{2}+n)}
2^{j(r_{2}-\sigma-s-s'+\frac{n}{p})}
||c||_{a^{s}(\dot{e}^{s'}_{pq})^{\sigma}_{x_{0}}}\times 
\\
&&
\sum_{l(P)=2^{-j}}(1+2^{j}|x_{R}-x_{P}|)^{-(L-\sigma)}
\\
&\leq& C2^{-i(n+\sigma+s+s'-\frac{n}{p})}
||c||_{a^{s}(\dot{e}^{s'}_{pq})^{\sigma}_{x_{0}}}
 < \infty
\end{eqnarray*}
if $r_{2} > \sigma + s +s'-\frac{n}{p}$ and $L > n+\sigma$.

Similarly we give the estimations of  $I_{1}$ using Lemma 2 and the maximal operator $M_{t}$ as in the proof of Lemma 1, 
\begin{eqnarray*}
\lefteqn{
|I_{1}| \leq C
\sum_{j \geq i}\sum_{l(P)=2^{-j}}|c(P)||\langle \phi_{P}\ ,\ \varphi_{R} \rangle|
}
\\
&\leq& 
C\sum_{j \geq i}\sum_{l(P)=2^{-j}}|c(P)
|2^{-jn}2^{(i-j)r_{1}}(1+2^{i}|x_{R}-x_{P}|)^{-L}
\\
&\leq&
C\sum_{j \geq i}2^{-j(r_{1}+n)}2^{ir_{1}}
\sum_{l(P)=2^{-j}}|c(P)|(1+2^{i}|x_{R}-x_{P}|)^{-L}
\\
&\leq&
\sum_{j \geq i}
2^{-j(r_{1}+n-n/t+s')}2^{ir_{1}}2^{-in/t}
M_{t}(\sum_{l(P)=2^{-j}}2^{js'}|c(P)|\chi_{P})(x)\ \ \ \ \ \ \ \ \ \ \ \ \ \ \ \ \ \ \ \ \ \ \ \ (a)
\end{eqnarray*}
if $0 < t \leq 1, \ L > n/t$ and $x \in R$ with $l(R)= 2^{-i}$.

Using  monotonicity of $l^q$-norm when $0 < q < 1$, and 
 H${\rm \ddot o}$lder's inequality when $1 \leq q \leq \infty$,
we get the following result from the above (a), 
$$
|I_{1}|\leq C2^{-i(n+s')}\{\sum_{i \leq j}(M_{t}(\sum_{l(P)=2^{-j}}2^{js'}
|c(P)|\chi_{P})(x))^q \}^{1/q}
$$
if $r_{1}+n-n/t+s' > 0,\ 0 < q \leq \infty$ and $x \in R$.

Taking $L^p(R)$ norm and using the  Fefferman-Stein vector valued inequality, 
we have, for a dyadic cube $R$ with $l(R)=2^{-i}$ and $c \in 
a^{s}(\dot{f}^{s'}_{pq})^{\sigma}_{x_{0}}$,
\begin{eqnarray*}
\lefteqn{
|I_{1}|2^{-in/p} = ||I_{1}||_{L^p(R)} 
}
\\
&\leq& C 2^{-i(n+s')}
||\{\sum_{i \leq j}M_{t}(\sum_{l(P)=2^{-j}}2^{js'}|c(P)|\chi_{P})^q\}^{1/q}||
_{L^{p}(R)}
\\
&\leq& 
C2^{-i(n+s')}||\{\sum_{i \leq j}(\sum_{l(P)=2^{-j}}2^{js'}|c(P)|
\chi_{P})^q\}^{1/q}||_{L^{p}(R)}
\\
&\leq&
C2^{-i(n+s')}c(\dot{f}^{s'}_{pq})(R)
\\
&\leq& 
C 2^{-i(n+s')}(|x_{0}-x_{R}|+2^{-i})^{\sigma}2^{-is}
||c||_{a^{s}(\dot{f}^{s'}_{pq})^{\sigma}_{x_{0}}}
\\
&\leq&
C2^{-i(n+s'+s+\sigma)}||c||_{a^{s}(\dot{f}^{s'}_{pq})^{\sigma}_{x_{0}}} < \infty\end{eqnarray*}
if $0 < t < \min(p,q)$, $0 < p < \infty,\ 0 < q \leq \infty$.

Next we give the estimates for the B-type case. We will consider two cases.

Case 1: $1 \leq p < \infty$. 
Taking $L^p(R)$ norm of (a) and using the boundedness of the maximal operator 
$M_{t}$,
we have  for a dyadic cube $R$ with $l(R)=2^{-i}$ and $c \in 
a^{s}(\dot{b}^{s'}_{pq})^{\sigma}_{x_{0}}$,
\begin{eqnarray*}
\lefteqn{
|I_{1}|2^{-in/p} = ||I_{1}||_{L^p(R)}
}
\\
&\leq& C 2^{i(r_{1}-n/t)}
||\sum_{i \leq j}2^{-j(r_{1}+n-n/t+s')}M_{t}(\sum_{l(P)=2^{-j}}2^{js'}|c(P)|\chi_{P})||_{L^{p}(R)}
\\
&\leq& 
C2^{i(r_{1}-n/t)}\sum_{i \leq j}2^{-j(r_{1}+n-n/t+s')}
||M_{t}(\sum_{l(P)=2^{-j}}2^{js'}|c(P)|\chi_{P})||_{L^{p}(R)}
\\
&\leq& 
C2^{i(r_{1}-n/t)}\sum_{i \leq j}2^{-j(r_{1}+n-n/t+s')}
||\sum_{l(P)=2^{-j}}2^{js'}|c(P)|\chi_{P}||_{L^{p}(R)}
\\
&\leq& 
C2^{-i(n+s')}\{\sum_{i \leq j}
||\sum_{l(P)=2^{-j}}2^{js'}|c(P)|\chi_{P}||_{L^{p}(R)}^q\}^{1/q}
\\
&\leq&
C2^{-i(n+s')}c(\dot{b}^{s'}_{pq})(R)
\leq
C2^{-i(n+s'+s+\sigma)}||c||_{a^{s}(\dot{b}^{s'}_{pq})^{\sigma}_{x_{0}}} < \infty\end{eqnarray*}
where these inequalities follow if $r_{1}+n-n/t+s'> 0$ and $p>t$ using  monotonicity of $l^q$-norm when $0 < q < 1$ and 
 H${\rm \ddot o}$lder's inequality when $1 \leq q \leq \infty$. 

Case 2: $0 < p <1$. Similarly, 
we have  for a dyadic cube $R$ with $l(R)=2^{-i}$ and $c \in 
a^{s}(\dot{b}^{s'}_{pq})^{\sigma}_{x_{0}}$,
\begin{eqnarray*}
\lefteqn{
|I_{1}|2^{-in/p} = ||I_{1}||_{L^p(R)}
}
\\
&\leq& C 2^{i(r_{1}-n/t)}
||\sum_{i \leq j}2^{-j(r_{1}+n-n/t+s')}M_{t}(\sum_{l(P)=2^{-j}}2^{js'}|c(P)|\chi_{P})||_{L^{p}(R)}
\\
&\leq& 
C2^{i(r_{1}-n/t)}\times
\\
&&
\{\sum_{i \leq j}2^{-j(r_{1}+n-n/t+s')p}
||M_{t}(\sum_{l(P)=2^{-j}}2^{js'}|c(P)|\chi_{P})||_{L^{p}(R)}^p\}^{1/p}
\\
&\leq& 
C2^{i(r_{1}-n/t)}\{\sum_{i \leq j}2^{-j(r_{1}+n-n/t+s')p}
||\sum_{l(P)=2^{-j}}2^{js'}|c(P)|\chi_{P}||_{L^{p}(R)}^p\}^{1/p}
\\
&\leq& 
C2^{-i(n+s')}\{\sum_{i \leq j}
||\sum_{l(P)=2^{-j}}2^{js'}|c(P)|\chi_{P}||_{L^{p}(R)}^q\}^{1/q}
\\
&\leq&
C2^{-i(n+s')}c(\dot{b}^{s'}_{pq})(R)
\leq
C2^{-i(n+s'+s+\sigma)}||c||_{a^{s}(\dot{b}^{s'}_{pq})^{\sigma}_{x_{0}}} < \infty\end{eqnarray*}
where these inequalities follow if $r_{1}+n-n/t+s'> 0$ and $p>t$ using  monotonicity of $l^q$-norm when $ q < p$ and 
 H${\rm \ddot o}$lder's inequality when $p \leq q $. 

We also get the same estimate for the case $p=\infty$.
Hence we obtain that
for a dyadic cube $R$ and a sequence $c \in 
a^s(\dot{e}^{s'}_{pq})^{\sigma}_{x_{0}}$,  
$\sum_{\mathcal{D} \ni P}c(P)\langle \phi_{P}\ ,\ \varphi_{R} \rangle$
 is convergent 
if  $r_{1} > J-n-s',\ \ r_{2} > \sigma+s+s'-\frac{n}{p}$ and $L>\max(J,\ n+\sigma)$.

(ii)\ \ Let $I_{0}$ and $I_{1}$ be as in the proof of (i) and $R$ a dyadic cube  with $l(R)=2^{-i}$. We put $w_{j}(P)=(2^{-j}+|x_{P}-x_{0}|)^{-\sigma},\ \ \sigma \geq 0$. 
Note  that 
$$
|c(P)| \leq Cl(P)^{s+s'-n/p}w_{j}(P)^{-1}
||c||_{a^{s}(\tilde{\dot{e}}^{s'}_{pq})^{\sigma}_{x_{0}}}$$ 
for $c \in a^{s}(\tilde{\dot{e}}^{s'}_{pq})^{\sigma}_{x_{0}}$. 
We have, by Lemma 2
\begin{eqnarray*}
\lefteqn{
|I_{0}| \leq C
\sum_{i \geq j}\sum_{l(P)=2^{-j}}|c(P)||\langle \phi_{P}\ ,\ \varphi_{R} \rangle|
}
\\
&\leq& 
C\sum_{i \geq j}\sum_{l(P)=2^{-j}}|c(P)
|2^{-in}2^{(j-i)r_{2}}(1+2^{j}|x_{R}-x_{P}|)^{-L}
\\
&\leq&
C\sum_{i \geq j}2^{-i(r_{2}+n)}2^{j(r_{2}-s-s'+\frac{n}{p})}
||c||_{a^{s}(\tilde{\dot{e}}^{s'}_{pq})^{\sigma}_{x_{0}}}\times 
\\
&&
\sum_{l(P)=2^{-j}}w_{j}(P)^{-1}(1+2^{j}|x_{R}-x_{P}|)^{-L}
\end{eqnarray*}
\begin{eqnarray*}
&\leq&
C\sum_{i \geq j}2^{-i(r_{2}+n)}2^{j(r_{2}-s-s'+\frac{n}{p})}
||c||_{a^{s}(\tilde{\dot{e}}^{s'}_{pq})^{\sigma}_{x_{0}}}\times 
\\
&&
2^{-j\sigma}\sum_{l(P)=2^{-j}}(1+2^{j}|x_{R}-x_{P}|)^{-(L-\sigma)}
\\
&\leq&
 C2^{-i(n+\sigma+s+s'-\frac{n}{p})}
||c||_{a^{s}(\tilde{\dot{e}}^{s'}_{pq})^{\sigma}_{x_{0}}}
 < \infty
\end{eqnarray*}
if $r_{2} > \sigma + s +s'-\frac{n}{p}$ and $L > \sigma+n$.

For $I_{1}$, we have
\begin{eqnarray*}
\lefteqn{
|I_{1}| \leq C
\sum_{j \geq i}\sum_{l(P)=2^{-j}}|c(P)||\langle \phi_{P}\ ,\ \varphi_{R} \rangle|
}
\\
&\leq& 
C\sum_{j \geq i}\sum_{l(P)=2^{-j}}|c(P)
|2^{-jn}2^{(i-j)r_{1}}(1+2^{i}|x_{R}-x_{P}|)^{-L}
\\
&\leq& 
C\sum_{j \geq i}\sum_{l(P)=2^{-j}}|c(P)
|2^{-jn}2^{(i-j)r_{1}}w_{j}(P)w_{j}(P)^{-1}(1+2^{i}|x_{R}-x_{P}|)^{-L}
\\
&\leq&
C\sum_{j \geq i}2^{-j(r_{1}+n)}2^{ir_{1}}2^{-i\sigma}
\sum_{l(P)=2^{-j}}|c(P)|w_{j}(P)(1+2^{i}|x_{R}-x_{P}|)^{-(L-\sigma)}
\\
&\leq&
\sum_{j \geq i}
2^{-j(r_{1}+n-n/t+s')}2^{i(r_{1}+\sigma-n/t)}
M_{t}(\sum_{l(P)=2^{-j}}2^{js'}w_{j}(P)|c(P)|\chi_{P})(x).
\end{eqnarray*}
By the same way as in the proof of (i), we get
$$
|I_{1}|\leq Cc(\tilde{\dot{e}}^{s'}_{pq})^{\sigma}_{x_{0}}(R)\leq C||c||_{a^{s}(\tilde{\dot{e}}^{s'}_{pq})^{\sigma}_{x_{0}}} < \infty$$
if $r_{1} > J-n-s'$ and $L > \sigma+J$.
We also obtain the same estimate for the case $\sigma < 0$.
 \qed

\bigskip

\begin{center}
{\bf 4. Characterizations}
\end{center}

 \bigskip

Let $\phi$ be a Schwartz function satisfying (1.1) and (1.2).
 We set 

\noindent
$\phi_{j}(x)= 2^{jn}\phi(2^j x)$ 
and 
$\phi_{Q}(x)=
\phi(l(Q)^{-1}(x-x_{Q}))$ 
for each dyadic cube $Q$.

\bigskip

{\bf Remark 1.}\  (See [7:\ \ Lemma(6.9)]). \ \ 
Note that for a given $\phi$ satisfying (1.1) and (1.2) 
there exists a Schwartz function $\varphi$ satisfying 
the same conditions (1.1) and (1.2)  such that 

\bigskip

$
\sum_{j \in {\Bbb Z}}\hat{\varphi}_{j}(\xi)
\hat{\phi}_{j}(\xi)=1$ if $\xi \neq 0.
$
\bigskip

\noindent
Hence we have
the $\varphi$-transform [5; Lemma 2.1]
$$
f=\sum_{R}l(R)^{-n}\langle f ,\ \varphi_{R} \rangle 
\phi_{R}
$$
for $f \in \mathcal{S}'$.
\bigskip

{\bf Theorem\ 1.}\  \ {\it 
For $s, \ s',\ \sigma \in {\Bbb R}$,\  $0 < p,\ q \leq \infty$,\  $x_{0} \in \mathbb{R}^n$ and $\phi \in \mathcal{S}$ satisfying (1.1) and (1.2), 
we have 

(i)\ \ 
$$
A^s(\dot{E}^{s'}_{pq})^{\sigma}_{x_{0}} = 
\{ f= \sum_{Q}c(Q)\phi_{Q} : 
\ \ (c(Q)) \in a^s(\dot{e}^{s'}_{pq})^{\sigma}_{x_{0}} \}.
$$  

(ii)\ \ 
$$
 A^s(\tilde{\dot{E}}^{s'}_{pq})^{\sigma}_{x_{0}} = 
\{ f= \sum_{Q}c(Q)\phi_{Q} : 
\ \ (c(Q)) \in a^s(\tilde{\dot{e}}^{s'}_{pq})^{\sigma}_{x_{0}} \}
$$
}

\bigskip

{\bf Remark\ 2.} 
(1)\ \ We can prove  that $
D\equiv
\{ f= \sum_{Q}c(Q)\phi_{Q} : 
\ \ c \in a^s(\dot{e}^{s'}_{pq})^{\sigma}_{x_{0}} \}$ and $\tilde{D}\equiv
\{ f= \sum_{Q}c(Q)\phi_{Q} : 
\ \ c \in a^s(\tilde{\dot{e}}^{s'}_{pq})^{\sigma}_{x_{0}} \}$

are independent of the choice of  $\phi \in {\mathcal S}$ satisfying the conditions  (1.1) and (1.2).
Indeed, suppose $\phi^1$ and $\phi^2$ are Schwartz functions satisfying (1.1) and (1.2), and the spaces $D(\phi^1)$ and $D(\phi^2)$ 
are defined by using  $\phi^1$ and $\phi^2$ in the place of $\phi$ respectively. We consider 
the $\varphi$-transform
$$
\phi^1_{P}=\sum_{R}l(R)^{-n}\langle \phi^1_{P}\ ,\ \varphi^2_{R} \rangle 
\phi^2_{R}.
$$
Then for $D(\phi^1) \ni 
f=\sum_{P}c(P)\phi^1_{P},\ c \in a^s(\dot{e}^{s'}_{pq})^{\sigma}_{x_{0}}$, 
we have   
$$
f=\sum_{P}c(P)\phi^1_{P}=\sum_{R}Ac(R)\phi^2_{R}
$$
where $A=\{l(R)^{-n}\langle \phi^1_{P}\ ,\ \varphi^2_{R} \rangle \}_{RP}$.
From Lemma 1 and Lemma 2, we see that  for 
$c \in a^s(\dot{e}^{s'}_{pq})^{\sigma}_{x_{0}}$,
$Ac \in a^s(\dot{e}^{s'}_{pq})^{\sigma}_{x_{0}}$. 
This shows that
$D(\phi^1) \subset D(\phi^2)$.
By the same argument, we see that  
$D(\phi^2) \subset D(\phi^1)$. These imply that the space 
$D$ is independent of the choice of $\phi$. This is also true for $\tilde{D}$

(2)\ \ 
From  Lemma 3, we note that
 for each  sequence $c \in a^s(\dot{e}^{s'}_{pq})^{\sigma}_{x_{0}}$, or  $c \in a^s(\tilde{\dot{e}}^{s'}_{pq})^{\sigma}_{x_{0}}$,  $\sum_{Q}c(Q)\phi_{Q}$ is convergent 
in ${\mathcal S}'_{k_{0}}$ 
where   $k_{0}=\max \{[\sigma+s + s'- \frac{n}{p}], -1\}$ and $\mathcal{S}_{-1}=\mathcal{S}, \mathcal{S}'_{-1}=\mathcal{S}'$, and 
$[x]$ is the greatest integer not greater than $x$.
Hence 
  we may  replace  $\mathcal{S}'_{\infty}$ in Definition 1 and Definition 2 
by $\mathcal{S}'_{k_{0}}$ (cf. [17]).

\bigskip

{\it Proof of Theorem 1} :\ \ 
We put $D\equiv \{ f= \sum_{Q}c(Q)\phi_{Q} : 
\ \ c \in a^s(\dot{e}^{s'}_{pq})^{\sigma}_{x_{0}} \}$. 
In order to prove $D \subset A^s(\dot{E}^{s'}_{pq})^{\sigma}_{x_{0}}$ we claim 
for a dyadic cube $P$ with $l(P)=2^{-j}$,
$$
c(\dot{E}^{s'}_{pq})(P)\leq Cc(\dot{e}^{s'}_{pq})(P).
$$
Let $D\ni f= \sum_{P}c(P)\phi_{P} :\ (c(P)) \in a^s(\dot{e}^{s'}_{pq})^{\sigma}_{x_{0}}$. 
Notice that  since $\phi_{i}*\phi \in \mathcal{S}$,
\begin{eqnarray*}
\lefteqn{|\phi_{i} * f(x)| = 
|\sum_{P} c(P) \phi_{i}*\phi_{P} (x)|
}
\\
&=&|\sum_{j=i-1}^{i+1}\sum_{l(P)=2^{-j}}c(P)\phi_{i}*\phi_{P}(x)|
\\
&=&|\sum_{j=i-1}^{i+1}\sum_{l(P)=2^{-j}}c(P)(\phi_{i-j}*\phi)_{P}(x)|
\\
&\leq& 
C \sum_{j=i-1}^{i+1}\sum_{l(P)=2^{-j}}|c(P)|(1+2^j|x-x_{P}|)^{-L} 
\end{eqnarray*}
for a large enough $L$. 
Hence we have, using the maximal operator $M_{t}$ as in the proof of Lemma 1,
\begin{eqnarray*}
\lefteqn{
\{\sum_{i\geq j} (2^{is'}|\phi_{i}*f|)^{q}\}^{1/q} \leq 
C\{\sum_{i\geq j} (2^{is'}\sum_{l(R)=2^{-i}}|\phi_{i}*f|\chi_{R})^q\}^{1/q}
}
\\
&\leq& 
C\{ \sum_{i \geq j}(2^{is'}
\sum_{l(R)=2^{-i}}(\sum_{l(R')=2^{-i+1}}|c(R')|
(1+2^{i-1}|x-x_{R'}|)^{-L}
\\ 
&+&
\sum_{l(R')=2^{-i}}|c(R')|
(1+2^{i}|x-x_{R'}|)^{-L}
\\
&+&
\sum_{l(R')=2^{-i-1}}|c(R')|
(1+2^{i+1}|x-x_{R'}|)^{-L})\chi_{R})^q \}^{1/q}
\end{eqnarray*}
\begin{eqnarray*}
\lefteqn{\leq
C\{ \sum_{i \geq j}(2^{is'}\times
}
\\
&&
\sum_{l(R)=2^{-i}}(\sum_{l(R')=2^{-i+1}}|c(R')|
(1+2^{i-1}|x_{R}-x_{R'}|)^{-L})\chi_{R})^{q}\}^{1/q}
\\
&+&
C\{\sum_{i \geq j}(2^{is'}
\sum_{l(R)=2^{-i}}(\sum_{l(R')=2^{-i}}|c(R')|
(1+2^{i}|x_{R}-x_{R'}|)^{-L})\chi_{R})^{q}\}^{1/q}
\\
&+&
C\{\sum_{i \geq j}(2^{is'}\times
\\
&&
\sum_{l(R)=2^{-i}}(\sum_{l(R')=2^{-i-1}}|c(R')|
(1+2^{i+1}|x_{R}-x_{R'}|)^{-L})\chi_{R})^q \}^{1/q}
\\
&\leq&
C\{ \sum_{i \geq j}(
\sum_{l(R)=2^{-i}}M_{t}(\sum_{l(R')=2^{-i+1}}2^{is'}|c(R')|\chi_{R'}
)\chi_{R})^{q}\}^{1/q}
\\
&+&
C\{\sum_{i \geq j}(
\sum_{l(R)=2^{-i}}M_{t}(\sum_{l(R')=2^{-i}}2^{is'}|c(R')|
\chi_{R'})\chi_{R})^{q}\}^{1/q}
\\
&+&
C\{\sum_{i \geq j}(
\sum_{l(R)=2^{-i}}M_{t}(\sum_{l(R')=2^{-i-1}}2^{is'}|c(R')|
\chi_{R'})\chi_{R})^q \}^{1/q}
\end{eqnarray*}
if $0< t \leq 1$ and $L > n/t$. 
Taking $L^p(P)$-norm and using the  Fefferman-Stein vector valued inequality, 
we have
\begin{eqnarray*}
\lefteqn{
c(\dot{F}^{s'}_{pq})(P)=||\{\sum_{i\geq j} (2^{is'}|\phi_{i}*f|)^{q}\}^{1/q}||_{L^p(P)} 
}
\\
&\leq&
C||\{ \sum_{i \geq j}(
M_{t}(\sum_{l(R')=2^{-i+1}}2^{is'}|c(R')|\chi_{R'}
))^{q}\}^{1/q}||_{L^p(P)}
\\
&+&
C||\{\sum_{i \geq j}(
M_{t}(\sum_{l(R')=2^{-i}}2^{is'}|c(R')|
\chi_{R'}))^{q}\}^{1/q}||_{L^p(P)}
\\
&+&
C||\{\sum_{i \geq j}(
M_{t}(\sum_{l(R')=2^{-i-1}}2^{is'}|c(R')|
\chi_{R'}))^q \}^{1/q}||_{L^p(P)}
\\
&\leq&
C||\{ \sum_{i \geq j}(
\sum_{l(R')=2^{-i+1}}2^{is'}|c(R')|\chi_{R'}
)^{q}\}^{1/q}||_{L^p(P)}
\\
&+&
C||\{\sum_{i \geq j}(
\sum_{l(R')=2^{-i}}2^{is'}|c(R')|
\chi_{R'})^{q}\}^{1/q}||_{L^p(P)}
\\
&+&
C||\{\sum_{i \geq j}(
\sum_{l(R')=2^{-i-1}}2^{is'}|c(R')|
\chi_{R'})^q \}^{1/q}||_{L^p(P)}
\\
&\leq&
C||\{\sum_{i \geq j}(
\sum_{l(R')=2^{-i}}2^{is'}|c(R')|
\chi_{R'})^{q}\}^{1/q}||_{L^p(P)} = Cc(\dot{f}^{s'}_{pq})(P)
\end{eqnarray*}
if $0< t <\min(p,q)$ and $1 \leq p < \infty$. 
If $0< p< 1$, we have 
\begin{eqnarray*}
\lefteqn{
c(\dot{F}^{s'}_{pq})(P)=||\{\sum_{i\geq j} (2^{is'}|\phi_{i}*f|)^{q}\}^{1/q}||_{L^p(P)} 
}
\\
&\leq&
C\{||\{ \sum_{i \geq j}(
M_{t}(\sum_{l(R')=2^{-i+1}}2^{is'}|c(R')|\chi_{R'}
))^{q}\}^{1/q}||^p_{L^p(P)}
\\
&+&
C||\{\sum_{i \geq j}(
M_{t}(\sum_{l(R')=2^{-i}}2^{is'}|c(R')|
\chi_{R'}))^{q}\}^{1/q}||^p_{L^p(P)}
\\
&+&
C||\{\sum_{i \geq j}(
M_{t}(\sum_{l(R')=2^{-i-1}}2^{is'}|c(R')|
\chi_{R'}))^q \}^{1/q}||^p_{L^p(P)}\}^{1/p}
\\
&\leq&
C\{||\{ \sum_{i \geq j}(
\sum_{l(R')=2^{-i+1}}2^{is'}|c(R')|\chi_{R'}
)^{q}\}^{1/q}||^p_{L^p(P)}
\\
&+&
C||\{\sum_{i \geq j}(
\sum_{l(R')=2^{-i}}2^{is'}|c(R')|
\chi_{R'})^{q}\}^{1/q}||^p_{L^p(P)}
\\
&+&
C||\{\sum_{i \geq j}(
\sum_{l(R')=2^{-i-1}}2^{is'}|c(R')|
\chi_{R'})^q \}^{1/q}||^p_{L^p(P)}\}^{1/p}
\\
&\leq&
C||\{\sum_{i \geq j}(
\sum_{l(R')=2^{-i}}2^{is'}|c(R')|
\chi_{R'})^{q}\}^{1/q}||_{L^p(P)} = Cc(\dot{f}^{s'}_{pq})(P).
\end{eqnarray*}
Similarly for B-type case we have 
\begin{eqnarray*}
\lefteqn{
c(\dot{B}^{s'}_{pq})(P)=
\{\sum_{i\geq j} (2^{is'}||\phi_{i}*f||_{L^p(P)} )^{q}\}^{1/q}\leq 
}
\\
&\leq&
C\{\sum_{i \geq j}
||M_{t}(\sum_{l(R')=2^{-i+1}}2^{is'}|c(R')|
\chi_{R'})||^q_{L^p(P)} \}^{1/q}
\\
&+&
C\{\sum_{i \geq j}
||M_{t}(\sum_{l(R')=2^{-i}}2^{is'}|c(R')|
\chi_{R'})||^q_{L^p(P)} \}^{1/q}
\\
&+&
C\{\sum_{i \geq j}
||M_{t}(\sum_{l(R')=2^{-i-1}}2^{is'}|c(R')|
\chi_{R'})||^q_{L^p(P)}\}^{1/q}
\\
&\leq&
C\{\sum_{i \geq j}
||\sum_{l(R')=2^{-i+1}}2^{is'}|c(R')|
\chi_{R'}||^q_{L^p(P)}\}^{1/q}
\\
&+&
C\{\sum_{i \geq j}
||\sum_{l(R')=2^{-i}}2^{is'}|c(R')|
\chi_{R'}||^q_{L^p(P)}\}^{1/q}
\\
&+&
C\{\sum_{i \geq j}
||\sum_{l(R')=2^{-i-1}}2^{is'}|c(R')|
\chi_{R'}||^q_{L^p(P)}\}^{1/q}
\\
&\leq&
C\{\sum_{i \geq j}
||\sum_{l(R')=2^{-i}}2^{is'}|c(R')|
\chi_{R'}||^q_{L^q(P)}\}^{1/q} = Cc(\dot{b}^{s'}_{pq})(P)
\end{eqnarray*}
if $1 \leq p < \infty$ and $0< t < p$. 
For $0< p < 1$ case , we obtain the same result by the same argument as in 
the proof of the F-type case. Moreover, for  $p=\infty$ case,  we have the same result.
Thus we have 
$$
c(\dot{E}^{s'}_{pq})(P)\leq Cc(\dot{e}^{s'}_{pq})(P)
$$
if $0 < p,q \leq \infty$. This implies 
$D \subset 
A^s(\dot{E}^{s'}_{pq})^{\sigma}_{x_{0}}$.

Next we will show the converse argument that
$A^s(\dot{E}^{s'}_{pq})^{\sigma}_{x_{0}} \subset D$.
We consider the $\varphi$-transform 
$f= \sum_{P}c(f)(P)\phi_{P}$, 
$c(f)(P)= l(P)^{-n}\langle f\ ,\ \varphi_{P} \rangle$ 
where $\phi$ and  $\varphi$ as in Remark 1. 
We will show  that $c(f)(P) \in a^s(\dot{e}^{s'}_{pq})^{\sigma}_{x_{0}}$. 
More precisely, we will see  that 
$$
c(f)(\dot{e}^{s'}_{pq})(P) \leq Cc(\dot{E}^{s'}_{pq})(P)
$$ for a dyadic cube $P$
with $l(P)=2^{-j}$ where notation
$c(f)(\dot{e}^{s'}_{pq})(P)$ means replacing 
with the sequence $c(f)(P)$ in the place of the sequence $c(P)$ 
in the definition of $c(\dot{e}^{s'}_{pq})(P)$. 
We define the sequence $\sup(f)(P)$ by setting 
$$
\sup(f)(P)= 
\sup_{P \ni y}|\varphi_{j}*f(y)|
$$  and  for $\gamma \in \mathbb{N}_{0}$, the sequence $\inf_{\gamma}(f)(P)$ and 
$t_{\gamma}(P)$ 
by 
$$
\inf_{\gamma}(f)(P)= \max\{\inf _{R \ni y}|\varphi_{j}*f (y)|: R \subset P, l(R)=2^{-(\gamma+j)} \},$$
$$\ \ 
t_{\gamma}(P)=\inf_{P \ni y}|\varphi_{j-\gamma}*f(y)|
$$
respectively. For a sequence $c(P)$ with $l(P)=2^{-j}$, we define a sequence $c^*(P)$ by 
$$
c^*(P) = \sum_{l(R)=2^{-j}}|c(R)|(1+2^{j}|x_{P}-x_{R}|)^{-L}
$$
for a sufficiently large $L$. 

We have, from Lemma A (i) in Appendix
\begin{eqnarray*}
\lefteqn{|c(f)(P)|= l(P)^{-n}|\langle f\ ,\ \varphi_{P} \rangle| 
=|\varphi_{j} * f (x_{P})| 
}
\\
&\leq&
 \sup(f)(P) 
\leq \sup(f)^*(P) 
\approx \inf_{\gamma}(f)^*(P)
\end{eqnarray*}
for $\gamma$ a large enough. 
 Hence by Lemma A  (ii) and (iii) in Appendix, we have  
\begin{eqnarray*}
\lefteqn{|c(f)(\dot{f}^{s'}_{pq})(P)| \leq C\inf_{\gamma}(f)^*(\dot{f}^{s'}_{pq})(P) \leq C \inf_{\gamma}(f)(\dot{f}^{s'}_{pq})(P)
}
\\
&\leq&
C||\{ \sum_{i\geq j}(\sum_{l(R)=2^{-i}}2^{is'}\inf_{\gamma}(f)(R)
\chi_{R})^{q} \}^{1/q}||_{L^p(P)} 
\\
\lefteqn{\leq
C||\{ \sum_{i\geq j}(\sum_{l(R)=2^{-i}}2^{is'}2^{\gamma L}
\sum_{R' \subset R, l(R')=2^{-(\gamma+i)}}t_{\gamma}^*(R')\chi_{R'}\chi_{R})^{q} \}^{1/q}||_{L^p(P)} 
}
\\
&\leq&
C2^{\gamma L}||\{ \sum_{i\geq j+\gamma}(2^{is'}2^{-\gamma s'}\sum_{l(R')=2^{-i}}t_{\gamma}(R')\chi_{R'})^{q} \}^{1/q}||_{L^p(P)} 
\\
&\leq&
C2^{\gamma( L- s')}||\{ \sum_{i\geq j+\gamma}(2^{is'}\sum_{l(R')=2^{-i}}
|\varphi_{i-\gamma}*f(y)|\chi_{R'})^{q} \}^{1/q}||_{L^p(P)} 
\\
&\leq&
C2^{\gamma( L- s')}||\{ \sum_{i\geq j}(2^{is'}2^{s'\gamma}
\sum_{l(R')=2^{-(i+\gamma)}}
|\varphi_{i}*f(y)|\chi_{R'})^{q} \}^{1/q}||_{L^p(P)} 
\\
&\leq&
C2^{\gamma L}||\{ \sum_{i\geq j}(2^{is'}
|\varphi_{i}*f(y)|)^{q} \}^{1/q}||_{L^p(P)}=Cc(\dot{F}^{s'}_{pq})(P). 
\end{eqnarray*}
For B-type case we can prove the same result by a similar argument as the above
$$
c(f)(\dot{b}^{s'}_{pq})(P) \leq C c(\dot{B}^{s'}_{pq})(P).
$$ 
Therefore we obtain 
$$
c(f)(\dot{e}^{s'}_{pq})(P) \leq C c(\dot{E}^{s'}_{pq})^{\varphi}(P)
$$
where the subscript $\varphi$  means replacing  $\phi$ with $\varphi$ in the definition of $c(\dot{E}^{s'}_{pq})(P)$. 
By Remark 2 (1),  
this implies 
$A^s(\dot{E}^{s'}_{pq})^{\sigma}_{x_{0}} \subset D$.
Hence we obtain $A^s(\dot{E}^{s'}_{pq}) = D$. We can prove (ii) by arguing as in the proof of (i).
\qed

\bigskip

{\bf Remark 3.} \ \ (1)\ \ From the above proof, we can see   that 

\bigskip

\noindent
$c(\dot{E}^{s'}_{pq})(P) \leq C c(\dot{e}^{s'}_{pq})(P)$  for 
$f=\sum_{P}c(P)\phi_{P},\ c \in a^s(\dot{e}^{s'}_{pq})^{\sigma}_{x_{0}}$, 

\noindent
$c(f)(\dot{e}^{s'}_{pq})(P) \leq C c(\dot{E}^{s'}_{pq})(P)$  for 
$c(f)(P)=l(P)^{-n}\langle f\ ,\ \phi_{P} \rangle,\ f \in A^s(\dot{E}^{s'}_{pq})^{\sigma}_{x_{0}}$ and, 

\bigskip

\noindent
$c(\tilde{\dot{E}}^{s'}_{pq})^{\sigma}_{x_{0}}(P) \leq 
C c(\tilde{\dot{e}}^{s'}_{pq})^{\sigma}_{x_{0}}(P)$  
for 
$f=\sum_{P}c(P)\phi_{P},\ c \in a^s(\tilde{\dot{e}}^{s'}_{pq})^{\sigma}_{x_{0}}$, 

\noindent
$c(f)(\tilde{\dot{e}}^{s'}_{pq})^{\sigma}_{x_{0}}(P) \leq C c(\tilde{\dot{E}}
^{s'}_{pq})^{\sigma}_{x_{0}}(P)$  for 
$c(f)(P)=l(P)^{-n}\langle f\ ,\ \phi_{P} \rangle,\ f \in A^s(\tilde{\dot{E}}^{s'}_{pq})^{\sigma}_{x_{0}}$.

\bigskip

\noindent
(2)\ \ $A^s(\dot{E}^{s'}_{pq})^{\sigma}_{x_{0}}$ and $A^s(\tilde{\dot{E}}^{s'}_{pq})^{\sigma}_{x_{0}}$ are independent of the choice of  $\phi \in {\mathcal S}$ satisfying the conditions  (1.1) and (1.2).

\bigskip

We have the following properties from Theorem 1. Hence we will investigate relations between our 2-microlocal spaces and the classical 2-microlocal spaces.

\bigskip

{\bf Proposition 1.}\ \ {\it 
 Suppose that 
\ $s, \ s,' \ \sigma \in {\mathbb R}$ and $x_{0} \in \mathbb{R}^n$.

\begin{enumerate}[{\rm (i)}]
\item \ When $\sigma < 0$, we have 

$A^s(\dot{E}^{s'}_{pq})^{\sigma}_{x_{0}}= \{ 0 \}$ for\ \ $0 < p, q \leq \infty$.
\item \ When $\sigma+s < 0$, 
we have 

$A^s(\dot{B}^{s'}_{pq})^{\sigma}_{x_{0}}= \{ 0 \}$  
for\ \ $0 < p,q \leq \infty$, 
$A^s(\dot{F}^{s'}_{pq})^{\sigma}_{x_{0}}= \{ 0 \}$ 
for\ \ $0 < p < \infty$,\  $0<q \leq \infty$.

When $\sigma+s+\frac{n}{q} < 0$, 
we have 

$A^s(\dot{F}^{s'}_{\infty q})^{\sigma}_{x_{0}}= \{ 0 \}$  
for\ \ $0 < q \leq \infty.
$
\item\ When $s < 0$, 
we have 

$A^s(\tilde{\dot{B}}^{s'}_{pq})^{\sigma}_{x_{0}}= \{ 0 \}$ 
for\ \ $0 < p,q \leq \infty$,
$A^s(\tilde{\dot{F}}^{s'}_{pq})^{\sigma}_{x_{0}}= \{ 0 \}$  
for\ \ $0 < p < \infty,\  0<q \leq \infty$. 

When $s+\frac{n}{q} < 0$,
 we have 

$A^s(\tilde{\dot{F}}^{s'}_{\infty q})^{\sigma}_{x_{0}}= \{ 0 \}$  for\ \ $0 < q \leq \infty$.
\end{enumerate}
}

\bigskip

{\it Proof} :\ \ These properties are shown from the fact that $c(\dot{e}^{s'}_{pq})(P)$,\ and 
 $c(\tilde{\dot{e}}^{s'}_{pq})^{\sigma}_{x_{0}}(P)$ 
are nondecreasing as $l(P)$ is increasing.
\qed

\bigskip

{\bf Proposition 2.}\ \ {\it 
 Suppose that 
\ $s, \ s,' \ \sigma \in {\mathbb R}$ and $x_{0} \in \mathbb{R}^n$.

\begin{enumerate}[{\rm (i)}]

\item \ We have

$A^{s}(\dot{B}^{s'}_{pq})^{\sigma} _{x_{0}}
\subset 
(\dot{B}^{s'}_{pq})^{s+\sigma}_{x_{0}}$ 
for\  $0 < p,q \leq \infty$,\

$A^{s}(\dot{F}^{s'}_{pq})^{\sigma} _{x_{0}}
\subset 
(\dot{F}^{s'}_{pq})^{s+\sigma}_{x_{0}}$ 
for\ \ $0 < p <\infty, 0<q \leq \infty$.

When $s \geq 0$, we have 

$A^s(\dot{F}^{s'}_{\infty q})^{\sigma}_{x_{0}}\subset (\dot{F}^{s'}_{\infty q})^{s+\sigma}_{x_{0}}$ 
 for\ \ $0< q \leq \infty$.
\item  
When $s \leq 0$, we have 

$A^s(\dot{B}^{s'}_{pq})^{\sigma}_{x_{0}}= (\dot{B}^{s'}_{pq})^{s+\sigma}_{x_{0}}$
for\  $0 < p,q \leq \infty$,\ 

$A^s(\dot{F}^{s'}_{pq})^{\sigma}_{x_{0}}= (\dot{F}^{s'}_{pq})^{s+\sigma}_{x_{0}}$
for\ \ $0 < p < \infty,\ 0<q \leq \infty$, 

$A^s(\dot{F}^{s'}_{\infty q})^{\sigma}_{x_{0}}\supset (\dot{F}^{s'}_{\infty q})^{s+\sigma}_{x_{0}}$ 
for\  \ $0< q \leq \infty$.

Particularly,
  when $\sigma \geq 0$ and $\sigma+s = 0$, 
we have 

$A^s(\dot{B}^{s'}_{pq})^{\sigma}_{x_{0}}= \dot{B}^{s'}_{pq}(\mathbb{R}^n)$ 
for\ \ $0 < p,q \leq \infty$, 

$A^s(\dot{F}^{s'}_{pq})^{\sigma}_{x_{0}}= \dot{F}^{s'}_{pq}(\mathbb{R}^n)$ 
for\ \ $0 < p <\infty,\ 0<q \leq \infty$. 
\item \ 
 When $\sigma < 0$, 
we have 
$A^s(\tilde{\dot{E}}^{s'}_{p q})^{\sigma}_{x_{0}}
\subset A^s(\dot{E}^{s'+\sigma}_{pq}) \subset 
\dot{E}^{s'+\sigma+s-\frac{n}{p}}_{\infty \infty}$ 
 for\ \ $0 < p,q \leq \infty$.
\end{enumerate}
}
{\it Proof} :
The properties (i) are followed from the fact that 
$$
l(Q)^{-\sigma}
\sup_{\mathcal{D} \ni P \subset 3Q}l(P)^{-s}c(\dot{e}^{s'}_{pq})(P)
\geq l(Q)^{-(\sigma+s)}c(\dot{e}^{s'}_{pq})(Q).
$$
\noindent
We can get the properties (ii) from the properties (i)
since, if $s \leq 0$, 
$$
a^s(\dot{e}^{s'}_{p q})^{\sigma}_{x_{0}}\supset (\dot{e}^{s'}_{p q})^{s+\sigma}_{x_{0}}. 
$$ 
The properties (iii) can be proved  since 
$$
l(P)^{-s}c(\tilde{\dot{e}}^{s'}_{pq})^{\sigma}_{x_{0}}(P)
 \geq l(P)^{-s}c(\dot{e}^{s'+\sigma}_{pq})(P) \geq 
l(P)^{-(s+s'+\sigma-\frac{n}{p})}|c(P)|.
$$
\qed

\bigskip

{\bf Proposition 3.}\ \ {\it 
 Suppose that 
\ $s, \ s,' \ \sigma \in {\mathbb R}$ and $x_{0} \in \mathbb{R}^n$. 

\begin{enumerate}[{\rm }]
\item\ 
When  $0 < q_{1} \leq q_{2} \leq \infty$, $0 <p \leq \infty$,
we have 

$A^{s}(\dot{B}^{s'}_{pq_{1}})^{\sigma} _{x_{0}}
\subset 
A^s(\dot{B}^{s'}_{pq_{2}})^{\sigma}_{x_{0}}$,
$A^{s}(\tilde{\dot{B}}^{s'}_{pq_{1}})^{\sigma} _{x_{0}}
\subset 
A^s(\tilde{\dot{B}}^{s'}_{pq_{2}})^{\sigma}_{x_{0}}$,

 and  
when  $0 < q_{1} \leq q_{2} \leq \infty$, $0 <p < \infty$,
 we have 

$A^{s}(\dot{F}^{s'}_{pq_{1}})^{\sigma} _{x_{0}}
\subset 
A^s(\dot{F}^{s'}_{pq_{2}})^{\sigma}_{x_{0}}$,
$A^{s}(\tilde{\dot{F}}^{s'}_{pq_{1}})^{\sigma} _{x_{0}}
\subset 
A^s(\tilde{\dot{F}}^{s'}_{pq_{2}})^{\sigma}_{x_{0}}$.

\end{enumerate}
}
{\it Proof} :\ \ These are consequences of the monotonicity of the $l^{p}$-norm.
\qed

\bigskip

{\bf Proposition 4.}\ \ {\it 
 Suppose that 
\ $s, \ s,' \ \sigma \in {\mathbb R}$ and $x_{0} \in \mathbb{R}^n$.

\begin{enumerate}[{\rm (i)}]
\item\ If 
$ 0 < p_{2} \leq p_{1} \leq \infty$ and $0 < q \leq \infty$, then 

$A^{s+\frac{n}{p_{1}}}(\dot{B}^{s'}_{p_{1}q})^{\sigma}_{x_{0}}
\subset 
A^{s+\frac{n}{p_{2}}}(\dot{B}^{s'}_{p_{2}q})^{\sigma}_{x_{0}}$, 
$A^{s+\frac{n}{p_{1}}}(\tilde{\dot{B}}^{s'}_{p_{1}q})^{\sigma}_{x_{0}}
\subset 
A^{s+\frac{n}{p_{2}}}(\tilde{\dot{B}}^{s'}_{p_{2}q})^{\sigma}_{x_{0}}$, 

 and,
if 
$ 0 < p_{2} \leq p_{1} < \infty$ and $0 < q \leq \infty$, then 

$A^{s+\frac{n}{p_{1}}}(\dot{F}^{s'}_{p_{1}q})^{\sigma}_{x_{0}}
\subset 
A^{s+\frac{n}{p_{2}}}(\dot{F}^{s'}_{p_{2}q})^{\sigma}_{x_{0}}$, 
$A^{s+\frac{n}{p_{1}}}(\tilde{\dot{F}}^{s'}_{p_{1}q})^{\sigma}_{x_{0}}
\subset 
A^{s+\frac{n}{p_{2}}}(\tilde{\dot{F}}^{s'}_{p_{2}q})^{\sigma}_{x_{0}}$, 

and, 
if
$ 0 < p \leq q \leq \infty$, then 

$A^{s}(\dot{F}^{s'}_{\infty q})^{\sigma}_{x_{0}}
\subset 
A^{s+\frac{n}{p}}(\dot{F}^{s'}_{pq})^{\sigma}_{x_{0}}$, 
$A^{s}(\tilde{\dot{F}}^{s'}_{\infty q})^{\sigma}_{x_{0}}
\subset 
A^{s+\frac{n}{p}}(\tilde{\dot{F}}^{s'}_{p q})^{\sigma}_{x_{0}}$.

\item\ 
If $0<  p,\ q \leq \infty$, then 

 $A^{s}(\dot{E}^{s'+\frac{n}{p}}_{pq})^{\sigma}_{x_{0}} 
 \subset A^{s}(\dot{E}^{s'}_{\infty\infty})^{\sigma}_{x_{0}} \subset 
 (\dot{E}^{s+s'}_{\infty\infty})^{\sigma}_{x_{0}}$, 

 $A^{s}(\tilde{\dot{E}}^{s'+\frac{n}{p}}_{pq})^{\sigma}_{x_{0}} 
 \subset A^{s}(\tilde{\dot{E}}^{s'}_{\infty\infty})^{\sigma}_{x_{0}} \subset 
 (\tilde{\dot{E}}^{s+s'}_{\infty\infty})^{\sigma}_{x_{0}}$.

\item\ 
If $s>0,\ 0<q< \infty,\ 0 < p \leq \infty$ 
or

  $s \geq 0,\ q=\infty, 
\ 0 < p \leq \infty$, 
then

$A^{s+\frac{n}{p}}(\dot{E}^{s'}_{pq})^{\sigma}_{x_{0}}=(\dot{E}^{s+s'}_{\infty\infty})^{\sigma}_{x_{0}}$, 
 and
$
A^{s+\frac{n}{p}}(\tilde{\dot{E}}^{s'}_{pq})^{\sigma}_{x_{0}}=(\tilde{\dot{E}}^{s+s'}_{\infty\infty})^{\sigma}_{x_{0}}$. 

\item\ If \ \ 
$ 0<  p_{1},\ p_{2},\ q \leq \infty$, then

$A^{\frac{n}{p_{1}}}(\dot{F}^{s'}_{p_{1}q})^{\sigma}_{x_{0}}
= 
A^{\frac{n}{p_{2}}}(\dot{F}^{s'}_{p_{2}q})^{\sigma}_{x_{0}}$,
$A^{\frac{n}{p_{1}}}(\tilde{\dot{F}}^{s'}_{p_{1}q})^{\sigma}_{x_{0}}
= 
A^{\frac{n}{p_{2}}}(\tilde{\dot{F}}^{s'}_{p_{2}q})^{\sigma}_{x_{0}}$.
\end{enumerate}
}
{\it Proof} :\ \ The properties (i) are  consequences of 
H$\ddot{\rm {o}}$lder's inequality. 
The properties (ii) are deduced from that 
$$
c(\dot{e}^{s'+\frac{n}{p}})(P)\geq l(R)^{-s'}|c(R)|
$$
for $R \subset P$.
 We will prove the properties  (iii). 
We see that 
$$
a^{s+\frac{n}{p}}(\dot{e}^{s'}_{pq})^{\sigma}_{x_{0}} \subset
(\dot{e}^{s'+s}_{\infty\infty})^{\sigma}_{x_{0}},
$$
since
$$
l(P)^{-(s+\frac{n}{p})}c(\dot{e}_{pq}^{s'})(P)\geq l(P)^{-(s'+s)}|c(P)|.
$$
In order to prove (iii), it is sufficient to prove  
$$
(\dot{e}^{s'+s}_{\infty\infty})^{\sigma}_{x_{0}} \subset 
a^{s+\frac{n}{p}}(\dot{e}^{s'}_{pq})^{\sigma}_{x_{0}}.
$$

\noindent
Since 
$$
c(\dot{e}^{s'}_{pq})(P) \leq C (\dot{e}^{s'+s}_{\infty\infty})(P)\times 
l(P)^{s+\frac{n}{p}}
$$ 

\noindent
if $s > 0$ and $0 < q < \infty$, hence 
we get the desired result. 
The properties (iv) is just [ 6: Corollary 5.7].  
\qed

\bigskip

{\bf Proposition 5.}\ \ {\it 
 Let $s,\ s' \in {\Bbb R},\ \ \sigma \geq 0,\ 0 < p,\ q \leq \infty\ and 
\ x_{0} \in \mathbb{R}^n$. We have

\begin{enumerate}[{\rm (i)}]
\item\ \ 
$
A^s(\dot{E}^{s'+\sigma}_{pq}) \subset 
A^s(\tilde{\dot{E}}^{s'}_{pq})^{\sigma}_{x_{0}} \subset 
 A^s(\dot{E}^{s'}_{pq})^{\sigma}_{x_{0}}$, 

\item \ \ $A^s(\tilde{\dot{B}}^{s'}_{p\ p\wedge q})^
{\sigma}_{x_{0}} \subset  A^s(\dot{F}^{s'}_{pq})^{\sigma}_{x_{0}}
 \subset  A^s(\dot{B}^{s'}_{p\ p\vee q})^{\sigma}_{x_{0}}$,

$A^s(\tilde{\dot{B}}^{s'}_{p\ p\wedge q})^{\sigma}_{x_{0}} \subset 
 A^s(\tilde{\dot{F}}^{s'}_{pq})^{\sigma}_{x_{0}} \subset 
  A^s(\dot{B}^{s'}_{p\ p\vee q })^{\sigma}_{x_{0}}$,

\item \ \ 
$
(\dot{E}^{s'}_{\infty \infty})^{\sigma}_{x_{0}}=
  (\tilde{\dot{E}}^{s'}_{\infty\infty})^{\sigma}_{x_{0}}.  
$  
Particularly, if $\frac{n}{p} < s$, 
then

$A^{s}(\dot{E}^{s'}_{pq})^{\sigma}_{x_{0}}=
A^{s}(\tilde{\dot{E}}^{s'}_{pq})^{\sigma}_{x_{0}}$.

\end{enumerate}
}
{\it Proof} :\ \ We obtain the properties (i) from the fact that  
$$
c(\tilde{\dot{e}}^{s'}_{pq})^{\sigma}_{x_{0}}(P) \leq Cc(\dot{e}^{\sigma+s'}_{pq})(P), 
$$
and 
$$
l(Q)^{-\sigma}l(P)^{-s}c(\dot{e}^{s'}_{pq})(P) \leq Cl(P)^{-s}c(\tilde{\dot{e}}^{s'}_{pq})^{\sigma}_{x_{0}}(P)
$$
since $l(Q)^{-\sigma} \leq (l(P)+|x_{0}-x_{P}|)^{-\sigma}$ for $P \subset 3Q$.

 The  embedding properties (ii)  can be proved   
from 

$c(\dot{f}^{s'}_{pq})(P) \leq c(\dot{b}^{s'}_{pq})(P)$\ \ 
for $0 < q \leq p \leq \infty$,

$c(\tilde{\dot{f}}^{s'}_{pq})(P) \leq c(\tilde{\dot{b}}^{s'}_{pq})(P)$\ \  
for $0 < q \leq p \leq \infty$,

$c(\dot{b}^{s'}_{pq})(P) \leq c(\dot{f}^{s'}_{pq})(P)$\ \ 
for $0< p \leq q \leq \infty$, 

$c(\tilde{\dot{b}}^{s'}_{pq})(P) \leq c(\tilde{\dot{f}}^{s'}_{pq})(P)$\ \  
for $0< p \leq q \leq \infty$ 

\noindent
by Minkowski's inequality, and the monotonicity (cf. Triebel[23: 2.3.2 Proposition 2]).

To prove the properties (iii), it is sufficient to see from properties (i),

$$
(\dot{e}^{s'}_{\infty \infty})^{\sigma}_{x_{0}} \subset
  (\tilde{\dot{e}}^{s'}_{\infty\infty})^{\sigma}_{x_{0}},  
$$  

\noindent 
 We consider  any dyadic cube $R$ with $l(R)=2^{-i}$ and 
 dyadic cubes $Q_{l}$ with $x_{0} \in Q_{l}$ and $l(Q_{l})=2^{-l},\ i \geq l$ such that 
 $Q_{i} \subset \cdots  \subset Q_{l} \subset Q_{l-1} \subset \cdots $ and 
 $\cup_{i \geq l}Q_{l}=\mathbb{R}^n$. We set 
 $Q_{l}^0 \equiv 3Q_{l} \setminus 3Q_{l+1},\ i >l$ 
 and $Q_{i}^0 \equiv 3Q_{i}$. 
 We divide into two cases. 
 
 Case (a): $R \subset Q_{l}^0,\ i > l$ case. Then we have 
 $2^{-i}+|x_{0}-x_{R}| \geq C 2^{-l},$
 
 Case (b):  $R \subset Q_{i}^0$. Then we have 
 $2^{-i}+ |x_{0} -x_{R}| \geq 2^{-i}$. 
 
\noindent  
In the case (a) we have
\begin{eqnarray*}
\lefteqn{
2^{is'}|c(R)|(2^{-i}+|x_{0}-x_{R}|)^{-\sigma} \leq C 2^{is'}2^{l\sigma}|c(R)| 
}
\\ 
&\leq&
C \sup_{x_{0}\in Q}2^{l\sigma}\sup_{R \subset 3Q}2^{is'}|c(R)| < \infty.
\end{eqnarray*}
In the case (b) we have 
\begin{eqnarray*}
\lefteqn{2^{is'}|c(R)|(2^{-i}+|x_{0}-x_{R}|)^{-\sigma}  \leq 
 C 2^{is'}2^{i\sigma}|c(R)|   
}
\\
&\leq&
C \sup_{x_{0} \in Q}2^{l\sigma}\sup_{R \subset 3Q}2^{is'}|c(R)|< \infty.
\end{eqnarray*}
These complete the proof.
\qed

\bigskip

{\bf Proposition 6.}\ \ {\it 
 Suppose that 
\ $s, \ s,' \ \sigma \in {\mathbb R}$,\ $0 <\epsilon$ and $x_{0} \in \mathbb{R}^n$. 

\noindent
\begin{enumerate}[{\rm (i)}]

\item\ 
 $A^{s}(\dot{B}^{s'+\epsilon}_{pq_{1}})^{\sigma-\epsilon}_{x_{0}} 
\subset 
A^{s}(\dot{B}^{s'}_{pq_{2}})^{\sigma}_{x_{0}}$ for\ \  $0<p\leq \infty,\ 0< q_{1}, \ q_{2} \leq \infty$, 

$A^{s}(\dot{F}^{s'+\epsilon}_{pq_{1}})^{\sigma-\epsilon}_{x_{0}} 
\subset 
A^{s}(\dot{F}^{s'}_{pq_{2}})^{\sigma}_{x_{0}}$ for\ \  $0<p< \infty,\ 0< q_{1}, \ q_{2} \leq \infty$, 
\item\ 
$A^{s+\epsilon}(\dot{E}^{s'}_{pq})^{\sigma-\epsilon}_{x_{0}} 
\subset 
A^{s}(\dot{E}^{s'}_{pq})^{\sigma}_{x_{0}}$ for\ \  $0< p, \ q \leq \infty$,

\item\ \ \ 
 $A^{s-\epsilon}(\dot{B}^{s'+\epsilon}_{pq_{1}})^{\sigma}_{x_{0}} 
\subset 
A^{s}(\dot{B}^{s'}_{pq_{2}})^{\sigma}_{x_{0}}$, and 
$A^{s-\epsilon}(\tilde{\dot{B}}^{s'+\epsilon}_{pq_{1}})^{\sigma}_{x_{0}} 
\subset 
A^{s}(\tilde{\dot{B}}^{s'}_{pq_{2}})^{\sigma}_{x_{0}}$ 
for\ \  $0< p,\ q_{1}, \ q_{2} \leq \infty$,

 $A^{s-\epsilon}(\dot{F}^{s'+\epsilon}_{pq_{1}})^{\sigma}_{x_{0}} 
\subset 
A^{s}(\dot{F}^{s'}_{pq_{2}})^{\sigma}_{x_{0}}$, 
 and 
$A^{s-\epsilon}(\tilde{\dot{F}}^{s'+\epsilon}_{pq_{1}})^{\sigma}_{x_{0}} 
\subset 
A^{s}(\tilde{\dot{F}}^{s'}_{pq_{2}})^{\sigma}_{x_{0}}$ 
for\ \  $0< p < \infty,\ 0< q_{1}, \ q_{2} \leq \infty$.

\end{enumerate}

}

\bigskip

{\it Proof}.\ \ (ii) is obvious.
 (i) and (iii) are consequences of  H$\ddot{\rm o}$lder's inequality and monotonicity of the $l^{p}$-norm.
\qed

\bigskip

We recall the definitions  of smooth atoms and molecules. 

\bigskip
 
 {\bf Definition 6.} \ \ Let $r_{1},\ r_{2} \in \mathbb{N}_{0}, L >n $.
  A family of functions $m = (m_{Q})$ indexed by dyadic cubes $Q$ 
is called a family of $(r_{1}, r_{2}, L)$-smooth molecules  if 

\bigskip

(3.1)\ \ 
 $
|m_{Q}(x)| \leq C(1+l(Q)^{-1}|x-x_{Q}|)^{-\max(L, L_{2})}
$ for some $L_{2} > n+r_{2}$,

(3.2)\ \ 
$|\partial^{\gamma}m_{Q}(x)| \leq Cl(Q)^{-|\gamma|}(1+l(Q)^{-1}|x-x_{Q}|)
^{-L}
$
 for\ \  $0 <|\gamma| \leq r_{1}$, and

(3.3)\ \ $\int_{{\mathbb R}^n}  x^{\gamma}m_{Q}(x) dx =0$ for\ \  $|\gamma| < r_{2}$,
 
\noindent
where (3.2) is void when $r_{1}=0$, and (3.3) is void when $r_{2}=0$.

\noindent
A family of functions $a = (a_{Q})$ indexed by dyadic cubes $Q$ 
is  called a family of  $(r_{1}, r_{2})$-smooth atoms 
if 

\bigskip

(3.4) \ \ 
supp $a_{Q} \subset 3Q$ for each dyadic cube $Q$,

(3.5) \ \ 
 $
|\partial^{\gamma}a_{Q}(x)| \leq Cl(Q)^{-|\gamma|}
$
 for\ \  $|\gamma| \leq r_{1}$, and

 (3.6)\ \ $\int_{{\mathbb R}^n}  x^{\gamma}a_{Q}(x) dx =0$ for\ \  $|\gamma| < r_{2}$,

\noindent
where  (3.6) is void when $r_{2}=0$.

\bigskip

{\bf Theorem 2.}\ \ {\it 
Let $s, \ s',\ \sigma \in {\Bbb R},$, $0 < p,\ q \leq \infty$\ and $x_{0} \in \mathbb{R}^n$.
Let $r_{1}, r_{2} \in \mathbb{N}_{0}$  
and $L > n$.

(i)\ \  We assume that  $r_{1}$, $r_{2}$ and $L$ satisfy 
 
\bigskip

 {\rm (4.1)}\ \ $r_{1} > \max (s', \sigma+s+s'-\frac{n}{p})$, 

 {\rm (4.2)}\ \ $r_{2} > J-n-s'$,

  {\rm (4.3)}\  \ $L > \max(J, n+\sigma)$

\noindent
where  $J$ as in Lemma 1.
 Then we have

\begin{eqnarray*}
\lefteqn{
A^s(\dot{E}^{s'}_{pq})^{\sigma}_{x_{0}}
=
\{ f =\sum_{Q}c(Q)m_{Q} : 
}
\\ 
&&
(r_{1},r_{2},L){\rm -smooth\ molecules}\ (m_{Q}),
\ \ (c(Q)) \in a^s(\dot{e}^{s'}_{pq})^{\sigma}_{x_{0}} \}
\end{eqnarray*}

$=\{f= \sum_{Q}c(Q)a_{Q} : 
 (r_{1},r_{2}){\rm  -smooth\  atoms}\ (a_{Q}),
\ \ (c(Q)) \in a^s(\dot{e}^{s'}_{pq})^{\sigma}_{x_{0}} \}$.

 \bigskip

(ii)\ \ We assume that  $r_{1}$, $r_{2}$ and $L$ satisfy 
 
\bigskip

 {\rm $(4.1)'$}\ \ $r_{1} > \max (s'+(\sigma\vee 0), (\sigma\vee 0)+s+s'-\frac{n}{p})$,

 {\rm $(4.2)'$}\ \ $r_{2} > J-n-s'-(\sigma\wedge 0)$,

 {\rm $(4.3)'$}\ \ $L > J+\sigma$

\noindent
where  $J$  as in Lemma 1.
Then we have

\begin{eqnarray*}
\lefteqn{
A^s(\tilde{\dot{E}}^{s'}_{pq})^{\sigma}_{x_{0}}
=
\{ f =\sum_{Q}c(Q)m_{Q} : 
}\\ 
&&
(r_{1},r_{2},L)-{\rm  smooth\ molecules}\ (m_{Q}),
\ \ (c(Q)) \in a^s(\tilde{\dot{e}}^{s'}_{pq})^{\sigma}_{x_{0}} \}
\\
&=&
\{f= \sum_{Q}c(Q)a_{Q} : 
\\ 
&&
(r_{1},r_{2})-{\rm  smooth\  atoms}\ (a_{Q}),
\ \ (c(Q)) \in a^s(\tilde{\dot{e}}^{s'}_{pq})^{\sigma}_{x_{0}} \},
\end{eqnarray*}

}

 \bigskip

{\bf Remark 4.}\ \    
From  Lemma 3, 
we remark that  
$f =\sum_{Q}c(Q)m_{Q}$ and 

\noindent
$f =\sum_{Q}c(Q)a_{Q}$ are convergent in 
 ${\mathcal S}'_{k_{0}}$.

\bigskip

{\it Proof of Theorem 2.}\ \  
(i)\ \ We put 

\noindent
$A\equiv \{f= \sum_{Q}c(Q)a_{Q} 
: \ \ (r_{1},r_{2}){\rm -smooth\  atoms}\ (a_{Q}),
 \ (c(Q)) \in a^s(\dot{e}^{s'}_{pq})^{\sigma}_{x_{0}} \}$, \\
$M \equiv \{ f =\sum_{Q}c(Q)m_{Q} : \ \ 
(r_{1},r_{2},L){\rm -smooth\ molecules}\ (m_{Q}),\\ 
 (c(Q)) \in a^s(\dot{e}^{s'}_{pq})^{\sigma}_{x_{0}} \}$. 

\noindent
Since an $(r_{1}, r_{2})$-atom is an $(r_{1}, r_{2}, L)$-molecule, it is easy to see that   $A \subset M$.
Let $M \ni f = \sum_{Q} c(Q) m_{Q}$ and  we consider the 
$\varphi$-transform 

$$
m_{Q}=\sum_{P} l(P)^{-n}\langle m_{Q}\ ,\ \varphi_{P} \rangle \phi_{P}
$$
where $\phi$ and $\varphi$ as in Remark 1.
Then we have

$$
f=\sum_{Q} c(Q) m_{Q} = \sum_{P} (Ac)(P) \phi_{P}
$$
where $A=\{ l(P)^{-n}\langle m_{Q}\ ,\ \varphi_{P} \rangle \}_{PQ}$.
Lemma 1 and Lemma 2 yield that 
$A$ is ($r_{1}, r_{2}+n, L$)-almost diagonal and 
$Ac \in a^s(\dot{e}^{s'}_{pq})^{\sigma}_{x_{0}}$
for $c \in a^s(\dot{e}^{s'}_{pq})^{\sigma}_{x_{0}}$.
Hence if we put $D \equiv \{ f= \sum_{Q}c(Q)\phi_{Q} : 
\ \ c \in a^s(\dot{e}^{s'}_{pq})^{\sigma}_{x_{0}} \}$, then  we see that 
$M \subset  D$. 

Using the argument similar to the proof of  
[6: Theorem 4.1](cf. [3: Theorem 5.9] 
or [4: Theorem 5.8]), 
for $D \ni f=\sum_{Q}c(Q)\phi_{Q},\ \  c \in 
a^s(\dot{e}^{s'}_{pq})^{\sigma}_{x_{0}}$.  
we see that there exist a family of $(r_{1}, r_{2})$-atoms $\{ a_{Q} \}$ 
and a sequence of coefficients $\{ c^{'}(Q) \} \in 
a^s(\dot{e}^{s'}_{pq})^{\sigma}_{x_{0}}$ such that 
$f=\sum_{Q}c(Q)\phi_{Q}=\sum_{Q}c^{'}(Q)a_{Q}$.
Hence we see that $D \subset A$. 
By Theorem 1,
 we obtain $A^s(\dot{E}^{s'}_{pq})^{\sigma}_{x_{0}} =D= M=A$.

(ii)\ \ We can prove by the same way in (i).

\qed

\bigskip

We recall the definition  of smooth wavelets. 

\bigskip
 
 {\bf Definition 7.}\ \ Let $r \in \mathbb{N}_{0}$ and $L >n$.
  A family of functions $\psi^{(i)}$ is called   a family of 
($r, L$)- smooth wavelets if
 $\{ 2^{nj/2}\psi^{(i)}(2^{j}x-k)\ 
(i = 1, \ldots, 2^n-1, \ j \in {\Bbb Z},\ k \in {\Bbb Z}^n) \}$ forms an 
orthonormal basis of $L^2({\Bbb R}^n)$ and  satisfies that for 
$\gamma \in {\mathbb N}_{0}^n$, 

\bigskip

(5.1)\ \ $|\psi^{(i)}(x)| \leq C(1+|x|)^{-\max(L, L_{0})}$ for some $ L_{0}> n+r$,

(5.2)\ \ $|\partial^{\gamma}\psi^{(i)}(x)| \leq C(1+|x|)^{-L}$ 
for $0 <|\gamma| \leq r$,

(5.3)\ \ $\int_{\mathbb{R}^n} \psi^{(i)}(x) x^{\gamma} dx =0$ for $|\gamma| <r$

\noindent
where (5.2) and (5.3) are void when $r=0$.

\noindent
We  denote 
$\psi_{Q}(x)=\psi(l(Q)^{-1}(x-x_{Q})), \ x_{Q}= 2^{-j}k$ for a dyadic cube 
$Q=[0, \ 2^{-j})^n + 2^{-j}k$. 
We will forget to write the index $i$ of the wavelet,  
which is of no consequence.

\bigskip

{\bf Theorem 3.}\ \  {\it
Let $s, \ s', \ \sigma \in {\Bbb R},\ x_{0} \in \mathbb{R}^n,\ 0 < p,\ q \leq \infty$ and $L > n$. 

(i)\ \ We assume that  a family of ($r, L$)- smooth wavelets $\psi$ satisfies

\bigskip

{\rm (6.1)}\ \  $r > \max(s',\ \sigma+s+s'-\frac{n}{p},\  J-n-s')$
 and

{\rm (6.2)}\ \ $L> \max(J,\ n+\sigma)$

\noindent
where $J$ is as in Lemma 1.

\bigskip

 Then we have

\bigskip

$A^s(\dot{E}^{s'}_{pq})^{\sigma}_{x_{0}}=
 \{ 
f = \sum_{Q}c(Q)\psi_{Q} : 
\ \ (c(Q)) \in a^s(\dot{e}^{s'}_{pq})^{\sigma}_{x_{0}} \}$.

\bigskip

(ii)\ \ We assume that a family of ($r, L$)- smooth wavelets $\psi$ satisfies

\bigskip

{\rm $(6.1)'$}\ \  $r > \max(s'+(\sigma\vee 0),(\sigma\vee 0)+s+s'-\frac{n}{p},\  J-n-s'-(\sigma\wedge 0))$
 and

{\rm $(6.2)'$}\ \ $L> J+\sigma $

\noindent

Then  
we have

\bigskip

$A^s(\tilde{\dot{E}}^{s'}_{pq})^{\sigma}_{x_{0}}=
 \{ 
f = \sum_{Q}c(Q)\psi_{Q} : 
\ \ (c(Q)) \in a^s(\tilde{\dot{e}}^{s'}_{pq})^{\sigma}_{x_{0}} \}$,

}

\bigskip

{\bf Remark 5.}\ \ 
(1)\ \ From Lemma 3 we note that $\sum_{Q}c(Q)\psi_{Q}$ is convergent in $\mathcal{S}_{k_{0}}'$. 

(2)\ \ 
we see that Theorem 3  
 is independent of the choice of  wavelets $\psi^{(i)}$ satisfying (6.1)-(6.2) 
or (6.1)'-(6.2)'  by Lemma 1 and Lemma 2.  

\bigskip

{\it Proof of Theorem 3} :\ \  
We put 
$W= \{ 
f = \sum_{Q}c(Q)\psi_{Q} : 
\ \ (c(Q)) \in a^s(\dot{e}^{s'}_{pq})^{\sigma}_{x_{0}} \}$.

\noindent
Let $W \ni f = \sum_{Q} c(Q) \psi_{Q}$ 
and  we consider the $\varphi$-transform

$$
\psi_{Q}=\sum_{P} l(P)^{-n}\langle \psi_{Q}\ ,\ \varphi_{P} \rangle \phi_{P}
$$

\noindent
where $\phi$ and $\varphi$ as in Remark 1.
Then we have

$$
f=\sum_{Q} c(Q) \psi_{Q} = \sum_{P} (Ac)(P) \phi_{P}\ \ \ \ \ \ \ \ \ \ \ \ \ \ \ \ \ \ \ \ \ \ \ \ \ \ \ \ \ \ \ \ \ \ \ \ \ \ \ \ \ \ \ \ \ \ \ \ (b)
$$
where $A=\{ l(P)^{-n}\langle \psi_{Q}\ ,\ \varphi_{P} \rangle \}_{PQ}$.
Lemma 1 and Lemma 2 yield that 
$A$ is ($r, r+n, L$)-almost diagonal and 
$Ac \in a^s(\dot{e}^{s'}_{pq})^{\sigma}_{x_{0}}$
for $c \in a^s(\dot{e}^{s'}_{pq})^{\sigma}_{x_{0}}$.
Hence 
by  Theorem 1, we see that 
$W \subset D =
 A^s(\dot{E}^{s'}_{pq})^{\sigma}_{x_{0}}$ where $D$ is as 
in the proof of Theorem 1.

Conversely, 
let $D \ni f = \sum_{Q} c(Q) \phi_{Q}$ and we consider the 
wavelet expansion

$$
\phi_{Q}=\sum_{P} l(P)^{-n}\langle \phi_{Q}\ ,\ \psi_{P} \rangle \psi_{P}.
$$

\noindent
Then we have

$$
f=\sum_{Q} c(Q) \phi_{Q} = \sum_{P} (Bc)(P) \psi_{P}
$$
where $B=\{ l(P)^{-n}\langle \phi_{Q}\ ,\ \psi_{P} \rangle \}_{PQ}$.
Lemma 1 and Lemma 2 yield that $B$ is 
\noindent
($r,r+n, L$)-almost diagonal and 
$Bc \in a^s(\dot{e}^{s'}_{pq})^{\sigma}_{x_{0}}$
for $c \in a^s(\dot{e}^{s'}_{pq})^{\sigma}_{x_{0}}$.
Hence by  Theorem 1, we see that 
$A^s(\dot{E}^{s'}_{pq})^{\sigma}_{x_{0}}=D \subset W$.  

We can prove (ii) by the same way in (i).
\qed

\bigskip

{\bf Remark 6.}\ \  
 For $f \in A^s(\dot{E}^{s'}_{pq})^{\sigma}_{x_{0}}$, the paring  
$\langle f\ ,\ \psi_{Q} \rangle$ is well-defined. More explicitly, we see  that for any $\{ \phi,\ \ \varphi \}$ as in Remark 1, by Lemma 3
$$
\langle f\ ,\ \psi_{Q} \rangle= 
\sum_{P}l(P)^{-n}\langle f\ ,\ \phi_{P} \rangle 
\langle \psi_{Q}\ ,\ \varphi_{P} \rangle=
\sum_{P}c(f)(P)\langle \psi_{Q}\ ,\ \varphi_{P} \rangle\ \ \ \ \ \ \ \ \ \ \ \ 
(c)
$$
\noindent
is convergent   sine  
$c(f)(P)=l(P)^{-n}\langle f\ ,\ \phi_{P} \rangle \in a^s(\dot{e}^{s'}_{pq})^{\sigma}_{x_{0}}$ (see  Remark 3 (1)).

Thus 
 for $f \in A^s(\dot{E}^{s'}_{pq})^{\sigma}_{x_{0}}$ 
we have a wavelet expansion $f= \sum_{Q}c(Q)\psi_{Q}$ in $\mathcal{S}_{k_{0}}'$  and 
the representation  $f= \sum_{Q}c(Q)\psi_{Q}$ is unique in $\mathcal{S}_{k_{0}}'$, 
that is,  
$c(Q)= l(Q)^{-n}\langle f\ ,\ \psi_{Q} \rangle$.
Hence we have 
$$
||f||_{A^s(\dot{E}^{s'}_{pq})^{\sigma}_{x_{0}}} \approx
||(c(Q))||_{a^s(\dot{e}^{s'}_{pq})^{\sigma}_{x_{0}}}.
$$

Indeed, from the above (c)  we see 
$$c(Q)= l(Q)^{-n}\langle f\ ,\ \psi_{Q} \rangle
=l(Q)^{-n}\sum_{P}c(f)(P)\langle \psi_{Q}\ ,\ \varphi_{P} \rangle=Ac(f)(Q)
$$ 
where $A=\{l(Q)^{-n}\langle \psi_{Q}\ ,\ \varphi_{P} \rangle\}_{QP}$. 
Therefore we have by Remark 3 (1) and Lemma 1, 
$$
||c||_{a^s(\dot{e}^{s'}_{pq})^{\sigma}_{x_{0}}} = ||Ac(f)||_{a^s(\dot{e}^{s'}_{pq})^{\sigma}_{x_{0}}} \leq C||c(f)||_{a^s(\dot{e}^{s'}_{pq})^{\sigma}_{x_{0}}}\leq 
C||f||_{A^s(\dot{E}^{s'}_{pq})^{\sigma}_{x_{0}}}.
$$ 
Conversely  from Remark 3 (1), Lemma 1 and (b) in the proof of Theorem 3  we have 
$$
||f||_{A^s(\dot{E}^{s'}_{pq})^{\sigma}_{x_{0}}} \leq 
||Ac||_{a^s(\dot{e}^{s'}_{pq})^{\sigma}_{x_{0}}}
\leq ||c||_{a^s(\dot{e}^{s'}_{pq})^{\sigma}_{x_{0}}}.
$$
Hence  
we have  an isomorphism 
$a^s(\dot{e}^{s'}_{pq})^{\sigma}_{x_{0}} \cong 
A^s(\dot{E}^{s'}_{pq})^{\sigma}_{x_{0}}$
in $\mathcal{S}_{k_{0}}'$
(not in $\mathcal{S}_{\infty}'$).

Similarly, 
for $f \in A^s(\tilde{\dot{E}}^{s'}_{pq})^{\sigma}_{x_{0}}$, 
the paring  
$\langle f\ ,\ \psi_{Q} \rangle$ is well-defined 
and we have an unique wavelet expansion $f= \sum_{Q}c(Q)\psi_{Q}$ in $\mathcal{S}_{k_{0}}'$. 
and we have 
$$||f||_{A^s(\tilde{\dot{E}}^{s'}_{pq})^{\sigma}_{x_{0}}} \approx
||(c(Q))||_{a^s(\tilde{\dot{e}}^{s'}_{pq})^{\sigma}_{x_{0}}}$$ 
(cf. [16: Theorems 4.1, 4.2 and Propositions 4.1, 4.2]).

\bigskip

\begin{center}
{\bf 5.  Applications}
\end{center}

\bigskip

{\bf Definition 8.}  
Let $\mathcal{T}$ be the space of Schwartz test functions 
($C^{\infty}$-functions with compact support)
and $\mathcal{T}'$ its dual.  
For an arbitrary $r_{1},\ r_{2} \in {\mathbb N}_{0}$ 
the Calder$\acute{{\rm o}}$n--Zygmund operator $T$ with an exponent $\epsilon 
> 0$  is 
a continuous linear operator $\mathcal{T} \rightarrow \mathcal{T}'$ such that 
its kernel $K$ off the diagonal $\{ (x,y) \in {\mathbb R}^n \times 
{\mathbb R}^n : x=y \}$ satisfies that

(7.1)\ \ $|\partial^{\gamma}_{1}K(x,y)| \leq C|x-y|^{-(n+|\gamma|)}$ 
for  $|\gamma| \leq  r_{1}$,

(7.2)\ \ $|K(x,\ y)-K(x,\ y')|
\leq  C|y-y'|^{r_{2} +\epsilon}|x-y|^{-(n+r_{2} + \epsilon)}$ 
if $2|y'-y| \leq |x-y|$,

(7.3)\ \ $|\partial^{\gamma}_{1}K(x,\ y)-\partial^{\gamma}_{1}K(x,\ y')|
\leq  C|y-y'|^{\epsilon}|x-y|^{-(n+|\gamma| + \epsilon)}$ 

if $2|y'-y| \leq |x-y|$ for $0 <|\gamma| \leq r_{1}$ 

(where this statement is void when $r_{1}=0$),

 $|\partial^{\gamma}_{1}K(x,\ y)-\partial^{\gamma}_{1}K(x',\ y)| 
\leq  C|x'-x|^{\epsilon}|x-y|^{-(n+|\gamma| + \epsilon)}$ 

if $2|x'-x| \leq |x-y|$ for $|\gamma| \leq r_{1}$,

\noindent
(where the subindex 1 stands for derivatives in the first variable)

(7.4)\ \ $T$ is bounded on $L^2({\Bbb R}^n)$.

\bigskip

We obtain the following theorem.

\bigskip

{\bf Theorem 4.} \ \ {\it 
 Let $s,\ s' , \ \sigma \in {\Bbb R},\   x_{0} \in \mathbb{R}^n$,
   $0< p,\ q \leq \infty$, $r_{1},\ r_{2} \in {\mathbb N}_{0}$ and 
 $J$ as in Lemma 1.  

(i)\ \ The Calder$\acute{o}$n--Zygmund operator $T$ with an exponent $\epsilon 
>\max(J,\ n+ \sigma)-n$  
 satisfying $T(x^{\gamma})=0$ for $|\gamma| \leq r_{1}$ and 
  $T^{*}(x^{\gamma})=0$ for $|\gamma| < r_{2}$,  
 is bounded on  
$A^s(\dot{E}^{s'}_{pq})^{\sigma}_{x_{0}}$
if $r_{1}$ and $r_{2}$ satisfy 

\bigskip

 {\rm $(8.1)$}\ \ $r_{1} > \max (s', \sigma+s+s'-\frac{n}{p})$,

 {\rm (8.2)}\ \ $r_{2} > J-n-s'$.

\bigskip

(ii)\ \ 
}The Calder$\acute{o}$n--Zygmund operator $T$ with an exponent $\epsilon 
>J -n+\sigma$  
 satisfying $T(x^{\gamma})=0$ for $|\gamma| \leq r_{1}$ and 
  $T^{*}(x^{\gamma})=0$ for $|\gamma| < r_{2}$,  
 is bounded on  
$A^s(\tilde{\dot{E}}^{s'}_{pq})^{\sigma}_{x_{0}}$ and 
$A^s(\tilde{E}^{s'}_{pq})^{\sigma}_{x_{0}}$ 
if $r_{1}$ and $r_{2}$ satisfy 

\bigskip

 {\rm $(8.1)'$}\ \ $r_{1} > \max (s'+(\sigma\vee 0), (\sigma\vee 0)+s+s'-\frac{n}{p})$,

 {\rm $(8.2)'$}\ \ $r_{2} > J-n-s'-(\sigma\wedge 0)$.

\bigskip

{\it Proof}.\ \  (i)\ 
The proof is similar to ones of [8].  
Let $f \in A^s(\dot{E}^{s'}_{pq})^{\sigma}_{x_{0}}$.
 We have a wavelet expansion 
$f = \sum_{Q}c(Q)\psi_{Q}: \  (c(Q)) \in a^s
(\dot{e}^{s'}_{pq})^{\sigma}_{x_{0}}$ from Theorem 3. 
We suppose that the wavelet $\psi$ is 
compactly supported with large enough smoothness 
by Remark 5 (2). Then there exists a positive constant c such that  
supp $\psi_{Q} \subset cQ$ for every dyadic cube $Q$.  

We claim that $Tf = \sum_{Q}c(Q)(T\psi_{Q}) \equiv \sum_{Q}c(Q)m_{Q}$ is convergent in 
$\mathcal{S}'_{k_{0}}$ and  
$||Tf||_{A^s(\dot{E}^{s'}_{pq})^{\sigma}_{x_{0}}} \leq  
C||f||_{A^s(\dot{E}^{s'}_{pq})^{\sigma}_{x_{0}}}$. 
Let $Q$ be a dyadic cube with $l(Q)=2^{-l}$.
From the assumption $T^{*}x^{\gamma}=0$ 
for $|\gamma| < r_{2}$ we have 
$\int_{{\mathbb R}^n}  x^{\gamma}m_{Q}(x) dx =0$ for\ \  $|\gamma| < r_{2}$.
  
 We choose a suitable large constant $C_{0}$. 
From Frazier--Torres--Weiss [8: Corollary 2.14],  when $|x-x_{Q}| < 2C_{0}2^{-l}$, we have 
$$
|\partial^{\gamma}m_{Q}(x)|\leq 
||\partial^{\gamma}m_{Q}||_{\infty} \leq 
C\sum_{|\alpha|\leq |\gamma|+1}2^{l(|\gamma|-|\alpha|)}2^{l|\alpha|}
||\partial^{\alpha}\psi_{Q}||_{\infty} $$
$$
\leq 
C 2^{l|\gamma|} \leq 
C l(Q)^{-|\gamma|}(1+l(Q)^{-1}|x-x_{Q}|)^{-L}
$$ 
for any $L \geq 0$ and $|\gamma| \leq r_{1}$.
 When $|x-x_{Q}| \geq 2C_{0}2^{-l}$, using the condition (7.2) in Definition 8, 
we obtain 
\begin{eqnarray*}
\lefteqn{
|m_{Q}(x)| = |\int_{\mathbb{R}^n}K(x,y)\psi_{Q}(y) dy|
}
\\
&=&
|\int_{\mathbb{R}^n}\bigl(K(x,y)-K(x,x_{Q})\bigr)\psi_{Q}(y) dy|  
\\
&\leq&
C\int_{|y-x_{Q}| \leq C_{0}2^{-l}}|K(x,y)-K(x,x_{Q})||\psi_{Q}(y)| dy
\\
&\leq& C \int_{|y-x_{Q}| \leq C_{0}2^{-l}}|y-x_{Q}|^{r_{2}+\epsilon}
|x-x_{Q}|^{-(n+r_{2}+\epsilon)} dy 
\\
&\leq& 
C(2^{l}|x-x_{Q}|)^{-(n+r_{2}+\epsilon)} 
\leq C(1+2^{l}|x-x_{Q}|)^{-(n+r_{2}+\epsilon)}.
\end{eqnarray*}
Moreover, using the condition (7.3) in Definition 8
for $0 < |\gamma| \leq r_{1}$, 
we have 
\begin{eqnarray*}
\lefteqn{
|\partial^{\gamma}m_{Q}(x)| \leq 
C\int_{|y-x_{Q}| \leq C_{0}2^{-l}}
|\partial^{\gamma}_{1}K(x,y)-\partial^{\gamma}_{1}K(x,x_{Q})||\psi_{Q}(y)| dy  
}
\\
&\leq&
 C\int_{|y-x_{Q}| \leq C_{0}2^{-l}}|y-x_{Q}|^{\epsilon}
|x-x_{Q}|^{-(n+|\gamma|+\epsilon)} dy 
\\
&\leq& 
C2^{-l(n+\epsilon)}|x-x_{Q}|^{-(n+|\gamma|+\epsilon)} 
\leq
C2^{l|\gamma|}(1+2^{l}|x-x_{Q}|)^{-(n+\epsilon)}.
\end{eqnarray*}
Therefore we observe that $m_{Q}=T\psi_{Q}$ is a molecule. More precisely 
$m_{Q}$ satisfies following properties: 

(3.1)\ \ 
 $
|m_{Q}(x)| \leq C(1+l(Q)^{-1}|x-x_{Q}|)^{-(n+r_{2}+\epsilon)}
$,

(3.2)\ \ 
$|\partial^{\gamma}m_{Q}(x)| \leq Cl(Q)^{-|\gamma|}(1+l(Q)^{-1}|x-x_{Q}|)
^{-(n+\epsilon)}
$
 for\ \  $0 <|\gamma| \leq r_{1}$, and

(3.3)\ \ $\int_{{\mathbb R}^n}  x^{\gamma}m_{Q}(x) dx =0$ for\ \  $|\gamma| < r_{2}$.

Hence by  Lemma 3, $Tf =\sum_{Q}c(Q)m_{Q}$ is convergent in 
$\mathcal{S}'_{k_{0}}$. For a wavelet expansion

$$
m_{Q}=\sum_{P} l(P)^{-n}\langle m_{Q}\ ,\ \psi_{P} \rangle \psi_{P},
$$
we have

$$
Tf=\sum_{Q} c(Q) m_{Q} = \sum_{P} (Ac)(P) \psi_{P}
$$
where $A=\{ l(P)^{-n}\langle m_{Q}\ ,\ \psi_{P} \rangle \}_{PQ}$.
From Lemma 1 and Lemma 2  the operator $A$ is bounded on $a^s(\dot{e}^{s'}_{pq})^{\sigma}_{x_{0}}$ if $r_{1}$ and $r_{2}$ satisfy (8.1) and (8.2) respectively. 
By Remark 6, 
it follows that 

$$||Tf||_{A^s(\dot{E}^{s'}_{pq})^{\sigma}_{x_{0}}} \approx 
||Ac||_{a^s(\dot{e}^{s'}_{pq})^{\sigma}_{x_{0}}}
\leq C||c||_{a^s(\dot{e}^{s'}_{pq})^{\sigma}_{x_{0}}} \approx 
C||f||_{A^s(\dot{E}^{s'}_{pq})^{\sigma}_{x_{0}}}.
$$ 
 This completes the proof of (i).

(ii)\ \  We can prove by the same way as in the proof of (i).
\qed

\bigskip

\bigskip

{\bf Definition 9.}

Let $\mu \in \mathbb{R}$. A smooth function $a$ defined on $\mathbb{R}^n \times \mathbb{R}^n$ is called to belong to the class $S^{\mu}_{1,1}(\mathbb{R}^n)$ if $a$ satisfies the following differential inequalities that for all $\alpha,\ \beta \in \mathbb{N}_{0}^n$,
$$
\sup_{x, \xi}(1+|\xi|)^{-\mu-|\alpha|+|\beta|}|\partial^{\alpha}_{x}
\partial^{\beta}_{\xi}a(x,\ \xi)|< \infty.
$$
$a(x, D)$ is the corresponding pseudo-differential operator such that 
$$
a(x, D)f(x)=\int_{\mathbb{R}^n}e^{ix\xi}a(x, \xi)\hat{f}(\xi)\ d\xi
$$
for $f \in \mathcal{S}$.

\bigskip

{\bf Theorem 5}. \ \ {\it 
 Let $s,\ s', \ \sigma \in {\Bbb R},\   x_{0} \in \mathbb{R}^n$,
   $0< p,\ q \leq \infty$.
Let $\mu \in \mathbb{R}$, $J$  as in Lemma 1 and $a \in S^{\mu}_{1,1}(\mathbb{R}^n)$.

(i)\ \ 
$a(x,D)$ is a continuous linear mapping from $A^s(\dot{E}^{s'}_{pq})^{\sigma}_{x_{0}}$ to $A^s(\dot{E}^{s'-\mu}_{pq})^{\sigma}_{x_{0}}$,  
 if $s' > J-n+\mu$  or if $a(x,\xi)=a(\xi)$. 

(ii)\ \ $a(x,D)$ is a continuous linear mapping  from $A^s(\tilde{\dot{E}}^{s'}_{pq})^{\sigma}_{x_{0}}$ to $A^s(\tilde{\dot{E}}^{s'-\mu}_{pq})^{\sigma}_{x_{0}}$  if $s' > J-n-(\sigma\wedge 0)+\mu$  or if $a(x,\xi)=a(\xi)$. 

}

\bigskip

{\it Proof}.\ \  We write $T\equiv a(x, D)$. Let $f \in 
A^s(\dot{E}^{s'}_{pq})^{\sigma}_{x_{0}}$. we consider a $\varphi$-transform 
$f = \sum_{P}c(P)\phi_{P}$  where $c(P)=c(f)(P)=l(P)^{-n}\langle f, \varphi_{P} \rangle$. 
We write  that $Tf= \sum_{P}c(P)m_{P}$, where $m_{P}=T\phi_{P}$. 
We see  

$$
m_{P}=\int e^{ix\xi}a(x,\xi)\hat{\phi}_{P}(\xi)\ d\xi.
$$
Then we have, using a change of variables,  
$$
m_{P}(x)=\int e^{i(x-x_{P})(2^j\xi)}a(x, 2^j\xi)\hat{\phi}(\xi)\ d\xi.
$$ 
By the fact that  $(1-\triangle_{\xi})^L(e^{ix\xi})=(1+|x|^2)^Le^{ix\xi}$ and an integration by parts, we obtain for $\gamma \in \mathbb{N}_{0}^n$
\begin{eqnarray*}
\lefteqn{\partial_{x}^{\gamma}m_{P}(x)
} 
\\
&=&
(1+(2^j|x-x_{P}|)^2)^{-L}\times
\\
&&
\int (1-\triangle_{\xi})^L
(e^{i2^j(x-x_{P})\xi})
\sum_{\delta \leq \gamma}(2^ji\xi)^{\delta}\partial_{x}^{\gamma-\delta}a(x, 2^j\xi)\hat{\phi}(\xi)\ d\xi 
\\
&=&
C(1+(2^j|x-x_{P}|)^2)^{-L}\times
\\
&&\int e^{i2^j(x-x_{P})\xi}
(1-\triangle_{\xi})^L
\sum_{\delta \leq \gamma}(2^ji\xi)^{\delta}\partial_{x}^{\gamma-\delta}a(x, 2^j\xi)
\hat{\phi}(\xi)\ d\xi 
\end{eqnarray*}
Thus we have 
\begin{eqnarray*}
\lefteqn{|\partial_{x}^{\gamma}m_{P}(x)|
} 
\\
&\leq&
C(1+2^j|x-x_{P}|))^{-2L}\times
\\
&&
\int \sum_{|\alpha+\beta+\tau|\leq 2L}\sum_{\delta \leq \gamma}2^{j|\delta|}2^{j|\beta|}
|\partial_{\xi}^{\alpha}(\xi)^{\delta}
||\partial_{\xi}^{\beta}\partial_{x}^{\gamma-\delta}a(x, 2^j\xi)||\partial_{\xi}^{\tau}\hat{\phi}(\xi)|\ d\xi 
\\
&\leq&
C(1+2^j|x-x_{P}|))^{-2L}\times
\\
&&
\int 
\sum_{|\alpha+\beta+\tau|\leq 2L}\sum_{\delta \leq \gamma}2^{j|\delta|}2^{j|\beta|}
|\partial_{\xi}^{\alpha}(\xi)^{\delta}|(1+2^j|\xi|)^{\mu+|\gamma|-|\delta|-|\beta|}|\partial_{\xi}^{\tau}\hat{\phi}(\xi)|\ d\xi 
\\
&\leq&
C2^{j\mu}2^{j|\gamma|}(1+2^j|x-x_{P}|))^{-2L}.
\end{eqnarray*}
Hence $m_{P}(x)$ satisfies
$$
|2^{-j\mu}\partial^{\gamma}m_{P}(x)| \leq C2^{j|\gamma|}(1+2^j|x-x_{P}|))^{-2L}.$$
for any $\gamma \in \mathbb{N}_{0}$ and any $L \geq 0$. 
We choose a suitable large $L$. 
For a wavelet transform
$$
2^{-j\mu}m_{P}=\sum_{R} l(R)^{-n}\langle 2^{-j\mu}m_{P}\ ,\ \psi_{R} \rangle \psi_{R},
$$
we have
$$
Tf=\sum_{P}2^{j\mu} c(P) (2^{-j\mu}m_{P}) = \sum_{R} A(2^{j\mu}c)(R) \psi_{R}
$$
where $A=\{ l(R)^{-n}\langle 2^{-j\mu}m_{P}\ ,\ \psi_{R} \rangle \}_{RP}$.
From Lemma 1 and Lemma 2, 
$A$ is bounded on $a^s(\dot{e}^{s'-\mu}_{pq})^{\sigma}_{x_{0}}$ if $s' > J-n+\mu$ or if $a(x,\xi)=a(\xi)$.  We remark that in the case $s' > J-n+\mu$, we do not assume the vanishing moment condition for $m_{P}$. 
But  in the case 
that $a(x,\xi)=a(\xi)$, we have the vanishing moment condition for $m_{P}$, indeed,   
$\int x^{\gamma}m_{P}(x)\ dx=C\partial^{\gamma}\hat{m_{P}}(0)=
C\partial^{\gamma}(\hat{\phi_{P}}\cdot a)(0)=0$ for any $\gamma \in \mathbb{N}_{0}$. 
From Remark 6 and Remark 3 (1), it follows that 

$||Tf||_{A^s(\dot{E}^{s'-\mu}_{pq})^{\sigma}_{x_{0}}} \leq 
C||A(2^{j\mu}c)||_{a^s(\dot{e}^{s'-\mu}_{pq})^{\sigma}_{x_{0}}}
\leq C||2^{j\mu}c||_{a^s(\dot{e}^{s'-\mu}_{pq})^{\sigma}_{x_{0}}}
\leq C||c||_{a^s(\dot{e}^{s'}_{pq})^{\sigma}_{x_{0}}}$
$\leq
C||f||_{A^s(\dot{E}^{s'}_{pq})^{\sigma}_{x_{0}}}.
$

Similarly we can prove for the other case (ii).  This completes the proof.
\qed

\bigskip

{\bf Corollary }.\ \ {\it
Let $s,\ s',\ \sigma\in {\Bbb R} ,\   x_{0} \in \mathbb{R}^n$,  $0< p,\ q \leq \infty$.

(i) \ \ 
Let  $\gamma \in \mathbb{N}^n_{0}$. 
Then the differential operator $\partial^{\gamma}$ 
 is  continuous  
 from $A^s(\dot{E}^{s'}_{pq})^{\sigma}_{x_{0}}$ to $A^s(\dot{E}^{s'-|\gamma|}_{pq})^{\sigma}_{x_{0}}$, 
 and from 
$A^s(\tilde{\dot{E}}^{s'}_{pq})^{\sigma}_{x_{0}}$ to $A^s(\tilde{\dot{E}}^{s'-|\gamma|}_{pq})^{\sigma}_{x_{0}}$. 

(ii)\ \  
Let $\mu \in {\Bbb R}$. 
The Bessel potential 
 $(1-\triangle)^{\mu/2}$ is continuous isomorphisms 
from $A^s(\dot{E}^{s'}_{pq})^{\sigma}_{x_{0}}$ onto $A^s(\dot{E}^{s'-\mu}_{pq})^{\sigma}_{x_{0}}$, and 
 from $A^s(\tilde{\dot{E}}^{s'}_{pq})^{\sigma}_{x_{0}}$ onto $A^s(\tilde{\dot{E}}^{s'-\mu}_{pq})^{\sigma}_{x_{0}}$ 
}

\bigskip

{\it Proof}.\ \  By Theorem 4,  since\ $\xi^{\gamma} \in S_{1,1}^{|\gamma|}$ and $(1+|\xi|^2)^{\mu/2} \in S_{1,1}^{\mu}$, the differential operator $\partial^{\gamma}$ is continuous mapping
from $A^s(\dot{E}^{s'}_{pq})^{\sigma}_{x_{0}}$ to $A^s(\dot{E}^{s'-|\gamma|}_{pq})^{\sigma}_{x_{0}}$,  and from 
$A^s(\tilde{\dot{E}}^{s'}_{pq})^{\sigma}_{x_{0}}$ to $A^s(\tilde{\dot{E}}^{s'-|\gamma|}_{pq})^{\sigma}_{x_{0}}$. 
and the Bessel potential are continuous mapping 
from $A^s(\dot{E}^{s'}_{pq})^{\sigma}_{x_{0}}$ to $A^s(\dot{E}^{s'-\mu}_{pq})^{\sigma}_{x_{0}}$,  and from 
$A^s(\tilde{\dot{E}}^{s'}_{pq})^{\sigma}_{x_{0}}$ to $A^s(\tilde{\dot{E}}^{s'-\mu}_{pq})^{\sigma}_{x_{0}}$. 
 To finish the proof of (ii) we need to show that the Bessel potential is  
surjective and one to one.   For 
$h \in  A^s(\dot{E}^{s'-\mu}_{pq})^{\sigma}_{x_{0}}$, we set $f = (1-\triangle)^{-\mu/2}h$. Then $h = (1-\triangle)^{\mu/2}f$. Similarly we can prove the 
other case. 
These complete the proof.
 \qed
\bigskip

\bigskip

\begin{center}
{\bf Appendix}
\end{center}

We will prove the Lemma used in the proof of Theorem 1.
For the notations see the proof of Theorem 1.
\bigskip

{\bf Lemma A}\ \  {\it We have,  for a dyadic cube $P$ with $l(P)=2^{-j}$,

\bigskip

{\rm (i)}\ 
$$
\sup(f)^*(P) \approx \inf_{\gamma}(f)^*(P)
$$
if $\gamma$ is sufficently large,

\bigskip

\rm (ii)\ 
$$
\inf_{\gamma}(f)(P)\chi_{P}\leq 2^{\gamma L}\sum_{R\subset P, l(R)
=2^{-(\gamma+j)}}t_{\gamma}^*(R)\chi_{R}
$$

\bigskip

{\rm (iii)}\ 
$$
c(\dot{e}^{s'}_{pq})(P) \approx c^*(\dot{e}^{s'}_{pq})(P), \ \ 
c(\tilde{\dot{e}}^{s'}_{pq})^{\sigma}_{x_{0}}(P) 
\approx c^*(\tilde{\dot{e}}^{s'}_{pq})^{\sigma}_{x_{0}}(P). 
$$
}

\bigskip

{\it Proof} : (i) is just [6, Lemma A.4]. 

(ii)\ \ Let $R_{0}$ and $R$  in $P$ be cubes with $l(R_{0})=l(R)=2^{-(\gamma+j)}$. It is sufficient to show 
$$
t_{\gamma}(R_{0}) \leq C2^{\gamma L}t_{\gamma}^*(R).
$$
Since 
$$
1 \leq C2^{\gamma L}(1+ 2^{\gamma +j}|x_{R}-x_{R_{0}}|)^{-L}, 
$$
we have 
\begin{eqnarray*}
\lefteqn{
t_{\gamma}(R_{0}) \leq Ct_{\gamma }(R_{0})2^{\gamma L}
(1+2^{\gamma +j}|x_{R}-x_{R_{0}}|)^{-L}
}
\\
&\leq& 
C2^{\gamma L}\sum_{l(R')=2^{-(\gamma+j)}}t_{\gamma}(R')
(1+2^{\gamma +j}|x_{R}-x_{R'}|)^{-L}=C2^{\gamma L}t_{\gamma}^*(R)
\end{eqnarray*}
(iii)\ \ It is sufficient to prove 
$$
c^*(\dot{e}^{s'}_{pq})(P) \leq C c(\dot{e}^{s'}_{pq})(P)
$$
since $|c(P)| \leq c^*(P)$.

Using the maximal operator $M_{t}\ (0< t \leq 1)$ as in the proof of Lemma 1 and the  Fefferman-Stein vector valued inequality, we have 
\begin{eqnarray*}
\lefteqn{
c^*(\dot{f}^{s'}_{pq})(P) = 
||\{ \sum_{i \geq j}
(2^{is'}\sum_{l(R)=2^{-i}}|c^*(R)|\chi_{R})^{q} \}^{1/q}||_{L^p(P)}
}
\\
&\leq& 
C||\{ \sum_{i \geq j}
(2^{is'}\times
\\
&&
\sum_{l(R)=2^{-i}}\sum_{l(R')=2^{-i}}|c(R')|
(1+2^{i}|x_{R}-x_{R'}|)^{-L}\chi_{R})^{q} \}^{1/q}||_{L^p(P)}
\\
&\leq& 
C||\{ \sum_{i \geq j}
(2^{is'}\sum_{l(R)=2^{-i}}M_{t}(\sum_{l(R')=2^{-i}}|c(R')|
\chi_{R'})\chi_{R})^{q} \}^{1/q}||_{L^p(P)}
\\
&\leq&
C||\{ \sum_{i \geq j}
(\sum_{l(R')=2^{-i}}2^{is'}|c(R')|\chi_{R'})^{q} \}^{1/q}||_{L^p(P)}
=Cc(\dot{f}^{s'}_{pq})(P)
\end{eqnarray*}
if $0< t < \min(p,q)$, $L > n/t$ and $0< p < \infty, 0< q \leq \infty$. 
For the  B-type case, we obtain the same result by the same argument as the above. 
Moreover, for  $p=\infty$ case,  we have the same result. We also have same result for $c(\tilde{\dot{e}}^{s'}_{pq})^{\sigma}_{x_{0}}(P)$ by same way as the above 

\bigskip

{\bf Ackowledgement}.\ \ 
The author would like to thank Prof. Yoshihiro Sawano for 
his encouragement and  many helpful remarks. The author thanks the referee for his/her valuable comments and his/her constructive suggestions.

\bigskip

{ Depatment of Mathematics}
 
{ Akita University,} 

{ 010-8502 Akita, Japan}

\bigskip

\bigskip

{\it e-mail:}  sakakoichi@gmail.com


\begin{thebibliography}
\rm
\bibitem{}
Ts. Batbold and Y. Sawano, 
{\em Decomposition for local Morrey spaces},
 {Eurasian\  Math.\  J.}\  {\bf 5} (2014), 9-45
\bibitem{Bony} 
J. M. Bony, 
{\em Second microlocalization and propagation of singularities for semilinear 
hyperbolic equations},  Hyperbolic equations and related topics, Academic Press (1988), 11-49.


\bibitem{Bownik} 
M. Bownik, 
{\em Atomic and molecular decompositions of anisotropic Besov spaces},  
{Math.\ Z.}\ {\bf 250} (2005), 539-571.

\bibitem{Bownik-Ho} 
M. Bownik and K-P Ho,
{\em Atomic and molecular decompositions of anisotropic Triebel--Lizorkin spaces},  {Trans.\ Amer.\ Math.\ Soc.}\ {\bf 358} (2005), 1469-1510.

\bibitem{Fraizer1} 
M. Frazier and B. Jawerth, 
{\em Decomposition of Besov spaces},  
{Indiana Univ.\ Math.\ J.}\ {\bf 34} (1985), 777-799.

\bibitem{Fraizer1} 
M. Frazier and B. Jawerth, 
{\em A discrete transform and decompositions of distribution spaces},  
{J.\ Func.\ Anal.}\ {\bf 93} (1990), 34-171.

\bibitem{Fraizer2} 
M. Frazier, B. Jawerth and G. Weiss,  
{\em Littlewood--Paley theory and the study of function spaces}, 
CBMS Regional Conf. Series in Math. {\bf 79}, Amer. Math. Soc., 1991.

\bibitem{Fraizer3} 
M. Frazier, R. Torres  and  G. Weiss, 
{\em The boundedness of Calder$\acute{\rm o}$n--Zygmund operators 
on the spaces 
$\dot{F}^{\alpha, q}_{p}$}, 
Revista Mat.\ Iberoamericana  {\bf 4}(1988), 41-70.


\bibitem{Jaffard1} 
S. Jaffard, 
{\em Pointwise regularity criteria}, 
{C.\ R.\ Acad.\ Sci.\ Paris}\   {\bf 1339} (2004), 757-762.


\bibitem{Jaffard2}
S. Jaffard, {\em Wavelet techniques for pointwise regularity}, 
{Annal.\ Fac.\ Sci.\ Toulouse}\  {\bf 15} (2006), 3-33.


\bibitem{Jaffard3}
S. Jaffard, {\em Pointwise regularity associated with function spaces and multifractal analysis},
Approximation and Probability, Banach  Center Publications, vol. 72, 
Institute of Mathematics, Polish Academy of Sciences (2006), 93-110.

\bibitem{Jaffard- Melot} 
S. Jaffard  and C. M${\rm \acute{e}}$lot, 
{\em Wavelet analysis of fractal boundaries Part 1 local exponents}, 
{Commun.\ Math.\ Phys.}\ {\bf 258} (2004), 513-539.


\bibitem{Kempka}
H. Kempka, {\em 2-microlocal Besov spaces and function spaces of 
varying smoothness}, 
{ Rev.\ Mat.\ Complut.}\ {\bf 22} (2009), 227-251.


\bibitem{Kempka2}
H. Kempka, {\em Atomic, molecular, wavelet decomposition of 
generalized 2-microlocal Besov spaces},
{ Funct.\ Approx.\ Comment.\ Math.}\ (2010), 171-208.



\bibitem{Komori}
Y. Komori-Furuya, K. Matsuoka, E. Nakai and Y. Sawano,  {\em Applications of Littlewood--Paley theory for 
$\dot{B}_{\sigma}$- Morrey spaces to the boundedness of integral operators}, 
{ J.\ Funct.\ spaces Appli.}\ (2013), 1-21.

\bibitem{Ky} 
G. Kyriazis, 
{\em Decomposition systems for function spaces}, 
 { Studia Math.}\ {\bf 157} (2003), 133-169.


\bibitem{Levy-Vehel Seuret}
J. Levy-Vehel and S. Seuret, 
{\em 2-microlocal formalism}, Fractal Geomery and Applications. 
A Jubilee of Benoit Mandelbrot, Proceedings of Symposia in Pure Mathematics, 
vol. 72, part 2 (2004), 153-215.

\bibitem{Meyer}
Y. Meyer, 
{\em Wavelet analysis, local Fourier analysis and 2-microlocalization}, 
{ Contemporary Math.}\ {\bf 189} (1995), 393-401.


\bibitem{Meyer}
Y. Meyer, 
{\em Wavelets, Vibrations and Scalings}, CRM Monograph Series. AMS vol. 9 (1998), 393-401.

\bibitem{Sawano-Tanaka}
Y. Sawano and H. Tanaka, 
{\em Decompositions of Besov-Morrey spaces and Triebel--Lizorkin-Morrey spaces},
{Math.\ Z.}\  {\bf 257} (2007), 871-905.


\bibitem{Sawano} 
Y. Sawano, D. Yang and W. Yuan, 
{\em New applications of Besov-type and Triebel--Lizorkin-type  spaces}, 
{J.\ Math.\ Anal.\ Appl.}\ {\bf 363} (2010), 73-85.

\bibitem{Sickel}
W. Sickel, D. Yang and W. Yuan, 
{\em Morrey and Campanato meet Besov, Lizorkin and Triebel}, 
Lect.\ Notes  in Math.\  vol. 2005, Springer, Berlin, 2010.


\bibitem{Tr1} 
H. Triebel,
{\em Theory of function spaces}, 
Birkh${\rm \ddot{a}}$euser, Basel, 1983.


\bibitem{Tr1} 
H. Triebel,
{\em Theory of function spaces II}, 
Birkh${\rm \ddot{a}}$euser, Basel, 1992.

\bibitem{Yang1}
D. Yang and W. Yuan, 
{\em New Besov-type spaces and Triebel--Lizorkin-type spaces}, 
{ Math.\ Z.} {\bf 265} (2010), 451-480.


\end{thebibliography}
 \end{document}